%% file: DeWit93.tex
\documentclass[fleqn,twoside,10pt,a4paper]{article}

\usepackage{amssymb}
\usepackage{epsfig}
\usepackage{supertabular}

\newtheorem{theorem}{Theorem}[section]
\newtheorem{definition}{Definition}[section]

\def\boxit#1{\vbox{\hrule\hbox{\vrule\kern3pt\vbox{\kern3pt#1\kern3pt%
                                           }\kern3pt\vrule}\hrule}}

\title{
  Gau{\ss} Cubature for the Surface of the Unit Sphere
}

\author{
  David~~De Wit \\
  \\
  BSc (Geology and Physics) 1987,
  BScAppHons (Geophysics) 1990 \\
  \\
  BSc (Mathematics) 1989,
  PGDipSc (Applied Mathematics) 1991 \\
  MScSt (Numerical Mathematics) 1992
}

\date{November 1993}

\begin{document}

\maketitle

\begin{abstract}
  Gau{\ss} cubature (multidimensional numerical integration) rules are
  the natural generalisation of the 1D Gau{\ss} rules. They are optimal
  in the sense that they exactly integrate polynomials of as high a
  degree as possible for a particular number of points (function
  evaluations). For smooth integrands, they are accurate,
  computationally efficient formulae.

  The construction of the points and weights of a Gau{\ss} rule
  requires the solution of a system of moment equations. In 1D, this
  system can be converted to a linear system, and a unique solution is
  obtained, for which the points lie within the region of integration,
  and the weights are all positive. These properties help ensure
  numerical stability, and we describe the rules as `good'. In the
  multidimensional case, the moment equations are nonlinear algebraic
  equations, and a solution is not guaranteed to even exist, let alone
  be good. The size and degree of the system grow with the degree of
  the desired cubature rule.  Analytic solution generally becomes
  impossible as the degree of the polynomial equations to be solved
  goes beyond $ 4 $, and numerical approximations are required. The
  uncertainty of the existence of solutions, coupled with the size and
  degree of the system makes the problem daunting for numerical
  methods.

  The construction of Gau{\ss} rules for (fully symmetric)
  $ n $-dimensional regions is easily specialised to the case of
  $ U_3 $, the unit sphere in 3D. Despite the problems described above,
  for degrees up to $ 17 $, good Gau{\ss} rules for $ U_3 $ have been
  constructed/discovered.
\end{abstract}

\raggedbottom

\pagebreak


\section{Introduction}

\subsection{Multidimensional Gau{\ss} Cubature}

Instead of directly considering the surface of the unit sphere $ U_3 $,
we will consider a more general case.  For $ n $-dimensional regions
$\mathcal{R}_n $, we will construct $ N $-point cubature rules
${\left\lbrace \mathbf{x}_i, w_i \right\rbrace}_{i=1}^N $ of the form:
\begin{eqnarray*}
  \int_{\mathcal{R}_n}
    \omega ( \mathbf{x} )
    f ( \mathbf{x} )
  d\mathbf{x}
  \;
  \approx
  \;
  \sum_{i=1}^N
    w_i
    f ( \mathbf{x}_i ).
\end{eqnarray*}
The $ \mathbf{x}_i \in \mathbb{R}^{n} $ are called the cubature
\emph{points}, and the $w_i \in \mathbb{R} $ are their respective
\emph{weights}.  We want rules that are as accurate as possible for a
given number of points (function evaluations). Smooth integrands may be
accurately approximated by polynomials, hence Gau{\ss} rules, which
exactly integrate all polynomials of as high a degree as possible, are
a natural choice.  For numerical stability, we also want rules with
positive weights, and points which lie within $\mathcal{R}_n $ (these
things are automatic in 1D). Up until about 1975, the best rules known
were the (Cartesian) product rules, which although good, are not
optimal.

The theory behind multidimensional Gau{\ss} rules for `fully symmetric'
$\mathcal{R}_n $ is a natural generalisation of the well-known 1D case.
It originates in a foundation paper by Mantel and Rabinowitz (1977)
\cite{MantelRabinowitz:77}.  The material has been applied to the case
of $ U_3 $ \cite{Keast:87}, although no computed rules have been
published.

We describe, in \S\ref{sec:GCforFS}, the theory behind the construction
of Gau{\ss} cubature rules for fully symmetric regions $\mathcal{R}_n $.
The material is specialised to the case of $ U_3 $ in
\S\ref{sec:U3Cubature}, which describes implementational details and
some computer programs. For degrees up to $ 17 $, good cubature rules
for $U_3 $ have been found, and these are listed in Appendix
\ref{app:GaussRulesU3}.  By comparison, the (non-optimal) product rules
for $ U_3 $ are simple to construct, and are not terribly inefficient.
We describe them, and some programs for their computation, in
\S\ref{sec:ProductRules}.  Alternative methods for cubature, such as
Monte Carlo and lattice methods are discussed in
\S\ref{sec:Alternatives}.  In general they are poor second choices when
Gau{\ss} rules are available.

Before progressing, we sketch an application of the use of cubature for
$U_3 $: the solution of the interior Dirichlet problem using a
boundary integral equation. More significant applications are from
statistics, where cubatures over many dimensions are required, and
efficiency is critical.


\subsection{Application: 3D Interior Dirichlet Problem}

\begin{definition}[3D Interior Dirichlet Problem]
  Given a smooth, bounded domain $ G \subset \mathbb{R}^{3} $, with
  boundary $ \partial G $, find $ u : G \to \mathbb{R} $ such that:
  \begin{enumerate}
  \item
    $ u $ satisfies Laplace's equation within $ G $, that is
    $ \nabla^2 u ( \mathbf{x} ) = 0 $,
    $ \forall \mathbf{x} \in G $.
  \item
    $ u $ is known on the boundary (the Dirichlet
    condition).  That is, there is a continuous function
    $ f : \partial G \to \mathbb{R} $ such that:
    $
      \forall \mathbf{x} \in \partial G
    $,
    $
      u ( \mathbf{x} ) = f ( \mathbf{x} )
    $.
  \end{enumerate}
\end{definition}

We will consider only domains $ G \subset \mathbb{R}^{3} $ of class
$C^2$ \cite[pp~21-22]{Kress:89}, which we will loosely call `smooth'.
(Their boundaries $ \partial G $ will be of class $ C^1 $.)  By
$C^k(G)$, we mean the set of $ k $ times continuously differentiable
real-valued functions defined on $ G $.  We will be interested in
functions $ u $ contained in $C ( \bar{G} ) \cap C^2 ( G ) $.

\pagebreak

The interior Dirichlet problem is a convenient model problem to work
with.  Its solution represents a potential function, that is readily
related to physical observables (e.g. electrostatic force).  It is
mathematically attractive, as the existence of a unique solution is
known.  Where the boundary $ \partial G $, and the boundary data $ f $
are simple, it may even be possible to find an analytic solution. In
general, however, this is not possible, and it is more sensible to
construct a numerical approximation.  Knowledge of the existence of a
unique solution, for even quite nonsmooth boundaries, greatly
encourages this.  Techniques to (approximately) solve the interior
Dirichlet problem may be able to be used as models for the solution of
more sophisticated boundary value PDEs.

One method of solving the interior Dirichlet problem is to reformulate
it in terms of a boundary integral equation on $ \partial G $.  This
reformulation is a natural choice: it is involved in a constructive
proof of the existence of a solution \cite{Kress:89}.  We begin by
defining the `fundamental solution of Laplace's equation in 3D' as the
function:%
\footnote{
  A different function is defined for a different number of dimensions,
  notably $ 2 $.
}
\begin{eqnarray*}
  \Phi ( \mathbf{x}, \mathbf{y} )
  =
  {\displaystyle\frac{1}{4 \pi | \mathbf{x} - \mathbf{y} |}}.
\end{eqnarray*}
For fixed $ \mathbf{y} \in \mathbb{R}^{3} $, $ \Phi ( \cdot, \mathbf{y} ) $ is harmonic
(satisfies Laplace's equation) in $ \mathbb{R}^{3} \setminus \left\lbrace \mathbf{y} \right\rbrace $.
Writing $ \mathbf{n} ( \mathbf{y} ) $ as the unit outward normal of $ \partial G $
at the point $ \mathbf{y} $, we construct the solution to the interior
Dirichlet problem in terms of the fundamental solution as follows:

\begin{theorem}[Solution to the interior Dirichlet problem]
  For $ \mathbf{x} \in G $, the \emph{double layer potential}
  \begin{eqnarray*}
    u ( \mathbf{x} )
    =
    \int_{\partial G}
      \phi ( \mathbf{y} )
      {\displaystyle\frac{\partial \Phi ( \mathbf{x}, \mathbf{y} )}{\partial \mathbf{n} ( \mathbf{y} )}}
    d \mathbf{s} ( \mathbf{y} ),
  \end{eqnarray*}
  with \emph{continuous density} $ \phi $ is a solution of the
  interior Dirichlet problem if $ \phi $ is the solution of the
  following integral equation, for $ \mathbf{x} \in \partial G $:
  \begin{eqnarray*}
    \phi ( \mathbf{x} )
    -
    2
    \int_{\partial G}
      \phi ( \mathbf{y} )
      {\displaystyle\frac{\partial \Phi ( \mathbf{x}, \mathbf{y} )}{\partial \mathbf{n} ( \mathbf{y} )}}
    d \mathbf{s} ( \mathbf{y} )
    =
    - 2
    f ( \mathbf{x} ).
  \end{eqnarray*}
\end{theorem}

A numerical approximation to the solution $ \phi $ of the BIE can be
used to construct a numerical approximation to the double layer
potential $ u $, and hence the solution to the Dirichlet problem
\cite{Atkinson:82b}.  Initially, we construct a cubature rule for
approximating integrals over $ \partial G $. This cubature rule is used
with a Galerkin technique to approximate the solution of the BIE for $
\phi $.  Lastly, $ u $ is approximated within $ G $ by application of
the cubature rule to the Galerkin approximation.  A cubature rule for
the manifold $ \partial G $ may be constructed by pointwise projection
from one for $ U_3 $. An advantage of this is that for a particular
number of desired points in the rule, a single rule for $ U_3 $ will
suffice for any manifold. A disadvantage is that although the rule may
be appropriate for $ U_3 $, the mapping process may reduce its
efficacy.  An even point density on $ U_3 $ may lose its regularity
when mapped, especially if $ \partial G $ is not concave. To illustrate
this, \cite{Atkinson:82a} uses the following region in experiments:
\begin{eqnarray*}
  \mathbf{x}
  =
  \sqrt
  {
    \cos ( 2 \theta )
    +
    \sqrt{c - \sin^2 ( 2 \theta )}
  }
  \left(
    \begin{array}{c}
      a \cos ( \phi ) \sin ( \theta ) \\
      b \sin ( \phi ) \sin ( \theta ) \\
      \cos ( \theta )
    \end{array}
  \right).
\end{eqnarray*}
A typical choice is $ ( a, b, c ) = ( 1, 2, 1.1 ) $.  As $ c $
decreases towards $ 1 $, the shape becomes less convex, eventually
becoming like a peanut, and numerical methods lose accuracy.

\vfill

\pagebreak


\section{Gau{\ss} Cubature for Fully Symmetric Regions}
\label{sec:GCforFS}

\subsection{Motivation for Cubature}

Let $\mathcal{R}_n $ be an $ n $-dimensional region contained in
$\mathbb{R}^{n} $; and $ \omega, f :\mathcal{R}_n \to \mathbb{R} $.  Consider the
numerical approximation of multiple integrals of the form:
\begin{eqnarray*}
  \int_{\mathcal{R}_n}
    \omega ( \mathbf{x} )
    f ( \mathbf{x} )
  d \mathbf{x}.
\end{eqnarray*}
Here, $ \omega $ is a weighting function, typically defined to contain the
`singularity' of the integrand. For example if $\mathcal{R}_n $ is
$\mathbb{R}^{3} $, a typical weighting function is
$\omega ( \mathbf{x} ) = e^{- {| \mathbf{x} |}^2} $, and for a wide class of
functions $ f $ (e.g. polynomials), the integral will exist.  Commonly
$\omega $ will be unity, and $\mathcal{R}_n $ will be spherically
symmetric about the origin. For example, if $\mathcal{R}_n $ is $ U_3 $,
we might have:
\begin{eqnarray*}
    \int_{U_3}
      x_1 x_2 x_3^2
    \;
    dx_1 \, dx_2 \, dx_3
    =
    0.
\end{eqnarray*}
Accurate approximation of multidimensional integrals is in general a
computationally expensive task. The rapid growth in expense with $ n $
has been called the `curse of dimensionality' (see
\S\ref{sec:ProductRulesIntro}).  Fortunately, for many applications,
$\mathcal{R}_n $ has some symmetry (consider $ U_3 $), and this can
greatly simplify both the construction and the application of rules.
Here, we consider the case where $\mathcal{R}_n $ has some symmetry
(not an issue in 1D!), and the construction of Gau{\ss} rules, which
minimise computational expense for a desired accuracy.  The philosophy
is that if $\mathcal{R}_n $ and $ \omega $ have a certain symmetry, then
it will be natural to place equally-weighted cubature points within
$\mathbb{R}^{n} $ according to the same symmetry. This symmetry should
considerably reduce the computations required to construct the rule.

\enlargethispage{\baselineskip}


\subsection{Full Symmetry}
\label{sec:FullSymmetry}

One of the most natural symmetries to conceive for $ \mathbb{R}^{n} $ is
\emph{full symmetry} in Cartesian coordinates.  A set of points is fully
symmetric if any point can be reached from any other point by a series
of orthogonal rotations about coordinate axes, and reflections in
coordinate planes. Observe that in 3D, the set of vertices of the unit
octahedron is fully symmetric; indeed a complete fully symmetric set of
points is an orbit under the action of the group $ G_8^* $ of
symmetries of the octahedron.  This group has order $ 48 $, and
includes both rotations and reflections.  In particular $ U_3 $ is
composed of complete sets of fully symmetric points.

The use of symmetry for cubature rules dates back to Russian works in
the 1960s, primarily in papers by Sobolev
\cite
{%
  Sobolev:62c,%
  Sobolev:62a,%
  Sobolev:62b%
}.
Other symmetries for 3D can be defined in terms of the symmetry groups
of other Platonic solids.  The immediately attractive case is that of
the dodecahedron/icosahedron, but this is harder to visualise than that
of the octahedron (where the vertices all lie on the coordinate axes),
and leads to more difficult algebra to disentangle. In any case,
alternative symmetries have little application beyond 3D, whilst the
notion of full symmetry generalises perfectly.  The notion of full
symmetry dates back to Lyness (1965) \cite{Lyness:65a}, and the use of
octahedral symmetry for cubature on $ U_3 $ was first put on a clear
foundation by Lebedev (1976) \cite{Lebedev:76}.  Shortly after this,
the foundation paper of Mantel and Rabinowitz (1977)
\cite{MantelRabinowitz:77} generalised the notion to arbitrary fully
symmetric domains, implicitly using the octahedral symmetry, although
not acknowledging the intellectual heritage.  This presentation of the
theory closely follows both \cite{MantelRabinowitz:77} and a
complementary paper by Keast and Lyness (1979) \cite{KeastLyness:79}.
To begin, we formalise the notion of full symmetry.

\begin{definition}[Full Symmetry between Two Points]
  $ \mathbf{x} $ and $ \mathbf{y} $ are a pair of \emph{fully symmetric} (FS)
  points in $ \mathbb{R}^{n} $, denoted $ \mathbf{x} \sim \mathbf{y} $, if $ \mathbf{y} $
  can be reached from $ \mathbf{x} $ by permutations and/or sign changes
  of the entries of $ \mathbf{x} $.
\end{definition}

Observe that ``$ \sim $'' is an equivalence relation, and thus it
induces a partition of $ \mathbb{R}^{n} $ into equivalence classes.
Suppose that
point $ \mathbf{y} $ has $ r $ non-zero coordinates, of which $ p $ are
distinct, and let the $ j $th distinct non-zero coordinate appear
$l_j $ times ($ j = 1, \dots, p $), so that
$l_1 + \dots + l_p = r $.  Trivially,
$0 \leqslant p \leqslant r \leqslant n $.  For this $ \mathbf{y} $, a point
$\mathbf{x} \sim \mathbf{y} $ can then be found such that:
\begin{equation}
  \mathbf{x}
  =
  (
    \underbrace{x_1, x_1, \dots, x_1}_{l_1 \mbox{times}},
    \underbrace{x_2, x_2, \dots, x_2}_{l_2 \mbox{times}},
    \dots,
    \underbrace{x_p, x_p, \dots, x_p}_{l_p \mbox{times}},
    0, 0, \dots, 0
  ).
  \label{eq:generator}
\end{equation}
Here we have chosen
$0 < x_1 \leqslant x_2 \leqslant \dots \leqslant x_p $.  This
$\mathbf{x} $ is called the \emph{generator} of the equivalence class,
and each class can be expressed uniquely in terms of a generator.  The
number of distinct points in the equivalence class containing generator
$ \mathbf{x} $ can be shown to be:
\begin{eqnarray*}
  \frac
  {
    2^n n!
  }{
    ( n - r )!
    l_1 !
    l_2 !
    \dots
    l_p !
  }.
\end{eqnarray*}
The number of elements in each class thus varies between $ 1 $ (for the
equivalence class $ [ \mathbf{0} ] $) and $ 2^n n! $ (in general).  The
exact number will be important for our work.

\begin{definition}[Full Symmetry for a Set of Points]
  A set of points $\mathcal{R}_n \subseteq \mathbb{R}^{n} $ is called
  \emph{fully symmetric} (FS) if $ \mathbf{x} \in\mathcal{R}_n $ and
  $ \mathbf{y} \sim \mathbf{x} $ imply $ \mathbf{y} \in\mathcal{R}_n $, that is the set
  $\mathcal{R}_n $ contains only complete equivalence classes.
\end{definition}

Clearly $ \mathbb{R}^{n} $ is an FS set of points. Other important examples of
FS sets of points (domains of integration) commonly found in the
literature are the $ n $-dimensional unit hypercube
$C_n = {[ -1, 1 ]}^n $; the $ n $-dimensional unit sphere
\\
\noindent
$
  S_n
  =
  \left\lbrace
    \mathbf{x} \in \mathbb{R}^{n}
    \; | \;
    x_1^2 + \dots + x_n^2
    \leqslant
    1
  \right\rbrace
$;
and its surface $ U_n $. In \S\ref{sec:U3Cubature}, we will be interested
in $ U_3 $.

\begin{definition}[Full Symmetry for a Function]
  Given an FS domain $\mathcal{R}_n $, a function
  $ g :\mathcal{R}_n \to \mathbb{R} $ is said to be a
  \emph{fully symmetric} (FS) \emph{function} if:  for all
  $ \mathbf{x}, \mathbf{y} \in\mathcal{R}_n $, $ \mathbf{x} \sim \mathbf{y} $ means that
  $ g ( \mathbf{x} ) = g ( \mathbf{y} ) $.
\end{definition}

It is not really necessary that the function be real-valued for this
definition to make sense, but this will be sufficient.  We will only
consider integrals over FS regions involving FS weight functions, and
will be approximating them using cubature rules on FS sets of points.
Where $\mathcal{R}_n $ is an $ n $-dimensional FS region, consider an
integrand $ f :\mathcal{R}_n \to \mathbb{R} $, which is not necessarily fully
symmetric, and an FS weight function
$\omega :\mathcal{R}_n \to \mathbb{R}^{0+} $ which is positive over a set of
positive volume:
\begin{eqnarray*}
  I [ f ]
  \equiv
  \int_{\mathcal{R}_n}
    \omega ( \mathbf{x} )
    f ( \mathbf{x} )
  d\mathbf{x}.
\end{eqnarray*}
This is to be approximated by an $ N $-point cubature rule
${\left\lbrace {\mathbf{x}}_i, w_i \right\rbrace}_{i=1}^N $, such that
$
  {\mathbf{x}}_i = {( x_{i1}, \dots, x_{in} )}^{\top} \in \mathbb{R}^{n}
$
and $w_i \in \mathbb{R} $, for $ i = 1, \dots, N $. The rule is then:
\begin{equation}
  I^N [ f ]
  =
  \sum_{i=1}^N
    w_i
    f ( {\mathbf{x}}_i ).
  \label{eq:integrule}
\end{equation}
Recall that the \emph{degree} of a polynomial in $ n $ variables is the
maximum \emph{sum} of the exponents in any of its terms, not the maximum
exponent of any one variable appearing in its terms.  Thus
$3 x_1^3 x_2^4 + x_1^5 x_3 $ is a polynomial in $ 3 $ (or more!)
variables, that is of degree $ 7 $ (not $ 5 $), whilst $ x_1^5 $ and
$x_1^2 x_2^3 $ are monomials of degree $ 5 $.  The rule $ I^N $ is
\emph{exact} for a function $ f $, if $ I^N [ f ] = I [ f ] $. If it
is exact for all polynomials in $ n $ variables of degree up to and
including $ k $, then it is called an integration rule of \emph{degree
(of exactness)} $ k $.  Cubature rules of high degrees of exactness
should accurately integrate smooth functions (which can be accurately
approximated by polynomials), but this does not necessarily carry over
to non-smooth functions.

\begin{definition}[Fully Symmetric Integration Rule]
  A cubature rule $ I^N $ is\\
  called a \emph{fully symmetric integration rule} if the evaluation
  points form an FS set, and all points in an FS equivalence class
  within the rule have the same weight.
\end{definition}

Note that this does \emph{not} require $ f $ to be FS, only that the
set of evaluation points $ {\mathbf{x}}_i $ are an FS set (and that the
weights are constant over all members of the same equivalence class).
An FS integration rule is then completely specified by a set of
generators and their corresponding weights.  This can greatly simplify
computations by reducing the number of points in a rule.  The degree of
an FS integration rule can be related to properties of polynomial
integrands $ f $:

\begin{enumerate}
\item
  If $ f $ is a monomial containing an odd power of some
  coordinate variable, then $ I [ f ] = I^N [ f ] = 0 $.
\item
  If $ f $ is a monomial containing only even powers of
  variables, then $ I [ f ] $ and $ I^N [ f ] $ depend only
  on the exponents and not on the ordering of the variables.
\end{enumerate}
Thus, an FS integration rule which is exact for all monomials of degree
up to and including $ 2 m $ is actually a rule of degree $ 2 m + 1 $,
and it suffices that it be exact for all monomials of the form:
\begin{eqnarray*}
  x_1^{2 k_1}
  x_2^{2 k_2}
  \dots
  x_{\mu}^{2 k_{\mu}},
\end{eqnarray*}
where $ 0 \leqslant \mu \leqslant n $ and
$1 \leqslant k_i \leqslant k_j $ for
$i \leqslant j $ and $ k_1 + \dots + k_{\mu} \leqslant m $. (Set
the monomial to $ 1 $ if $ \mu = 0 $.)

This rule can be written, for an appropriate set of generators
$X \subseteq {\left\lbrace {\mathbf{x}}_j \right\rbrace}_{j=1}^N $,
with elements $ \mathbf{x} $, each of which has an equivalence class
$[ \mathbf{x} ] $ with elements $ \mathbf{y} $:
\begin{equation}
  I^N [ f ]
  =
  \sum_{\mathbf{x} \in X}
    w_{\mathbf{x}}
    \sum_{\mathbf{y} \in [ \mathbf{x} ]}
      f ( \mathbf{y} ).
  \label{eq:eval}
\end{equation}
That is, to evaluate $ I^N [ f ] $:
\begin{enumerate}
\item
  For each equivalence class $ [ \mathbf{x} ] $, sum the function values
  over all elements of the class.
\item
  Multiply the weight associated with each generator by its respective
  equivalence class sum, and sum these multiples.
\end{enumerate}

\pagebreak


\subsection{Gau{\ss} Cubature for Fully Symmetric Regions}

Cubature rules that have a minimal number of points for a specified
degree are called `minimal rules'.

\begin{definition}[Fully Symmetric Minimal Rule]
  An $ N $-point FS rule of degree $ 2 m + 1 $ over an FS set
  $\mathcal{R}_n $ is a \emph{fully symmetric minimal} (FSM) rule if
  no other FS rule of degree $ 2 m + 1 $ over $\mathcal{R}_n $
  exists with less than $ N $ evaluation points.
\end{definition}

Note that an FSM rule is not necessarily unique. In 1D, FSM rules
\emph{are} the unique Gau{\ss} rules, for which the theory is
well-known.  We must be careful not to confuse Gau{\ss} rules with
product rules (see \S\ref{sec:ProductRules}), where Gau{\ss}-Legendre
rules are used as basic rules in the construction of (Cartesian)
product multidimensional cubature rules. Sometimes these rules are
called Gau{\ss} product rules. Whilst they are Gau{\ss} rules, in the
sense that they exactly integrate polynomials of as high a degree as
possible (but only in one dimension), they are \emph{not} minimal.
We would like our rule to possess a couple of important properties:

\begin{definition}[Good Rule]
  An integration rule $ \left\lbrace \mathbf{x}_i, w_i \right\rbrace $ over $\mathcal{R}_n $ is a
  \emph{good} rule if its evaluation points $ \mathbf{x}_i $ lie within
  $\mathcal{R}_n $, and its weights $ w_i $ are positive.
\end{definition}

The first of these conditions is familiar from 1D Gau{\ss} quadrature,
whilst the second is new.  These properties seem natural, however there
is nothing in the assumptions for cubature that demands them.  In 1D,
na{\"\i}ve approaches to quadrature, such as the (equally-spaced)
Newton-Cotes family, do not preserve the positivity of weights for
larger $ N $, whilst both these properties are satisfied by 1D Gau{\ss}
rules.  In higher dimensions, Gau{\ss} rules are not necessarily good.
The properties of good rules, in particular the latter, assist
numerical stability.  Whilst a rule may theoretically exactly integrate
the relevant polynomials, in practice, summation of large terms of
alternating sign tends to reduce numerical precision.  In
\S\ref{sec:U3Cubature}, we consider the case where $\mathcal{R}_n $ is $ U_3
$, so the internal points condition will simply require points to lie
on the surface of the unit sphere.

We use the concept of full symmetry to define several classes of rules
for the cubature of integrals with FS weight functions where the
\emph{moments}%
\footnote{
  That is, the integrals of appropriate polynomials over $\mathcal{R}_n $.
}
exist (commonly $ \omega \equiv 1 $). A fully symmetric cubature rule
of degree $ 2m + 1 $ exactly integrates all polynomials of degree up to
and including $ 2m $.  A fully symmetric minimal (FSM) rule does so
using a minimal number of cubature points.  A \emph{fully symmetric
good} (FSG) rule does so where all weights are positive and all points
lie within $\mathcal{R}_n $. A \emph{fully symmetric minimal good} (FSMG)
rule is both minimal and good. A \emph{fully symmetric good minimal}
(FSGM) rule is a good rule that is minimal, that is, although FS rules
on fewer points may exist, none are good.

We seek firstly FSMG rules, and if there are none of these, FSGM rules,
which always exist. At worst, FSGM rules are \emph{product rules} (see
\S\ref{sec:ProductRules}), which (generally) require
${[ 2 ( 2m + 1 ) - 1 ]}^n = {( 4m + 1 )}^n $ points to be of
degree $ 2m + 1 $.  For low-dimensional applications (small $ n $),
this may not be terribly inefficient, e.g. product rules for $ U_3 $
actually require $ 2 {( m + 1 )}^2 $ points.  Apart from such
considerations, a rule that is minimal (or almost so) might be
\emph{almost} good in the sense that negative weights are very small, or
that points are only just outside $\mathcal{R}_n $, or that failing
these conditions, errors in integrating polynomials are minor. The
quest for good rules may be an arduous search through these almost good
rules.


\subsection{Conditions for Gau{\ss} Cubature in 3D}
\label{sec:GCin3D}

The procedure for constructing Gau{\ss} rules becomes more complicated
as the dimension $ n $ increases. Here, we fully develop the 3D case,
specialising this to the case of $ U_3 $ in \S\ref{sec:U3Cubature}.  We
will set up a system of (moment) equations to express the fact that an
FS rule for $\mathcal{R}_3 $ will exactly integrate all polynomials of up
to a specified degree, \emph{without} presupposing the number or
distribution (beyond being an FS set) of points, or the sign of
weights.

To begin, we group generators into types, depending on the number of
zeros and repeated elements in the entries of their equivalence
classes.  There will always be the (not very interesting) class
$[ \mathbf{0} ] $, of unit multiplicity, and in general $ e $ other
types of classes, for a total of $ e + 1 $ types of classes. For small
$n $, this $ e $ is generally small, and may be found as the solution
to the following problem:%
\footnote
{
  The framing of this problem is the basis for the higher-dimensional
  analysis in \protect\cite{KeastLyness:79}.
}

Given a positive integer $ n $, $ e + 1 $ is the total number of
strings of length up to $ n $, with $ p $ distinct entries taken
from the positive integers $ 1, \dots, p $, of the form:
\begin{eqnarray*}
    \underbrace{1, 1, \dots, 1}_{l_1 \mbox{times}}, \,
    \underbrace{2, 2, \dots, 2}_{l_2 \mbox{times}}, \,
    \dots, \,
    \underbrace{p, p, \dots, p}_{l_p \mbox{times}}
    \qquad \qquad
    p \leqslant n,
\end{eqnarray*}
such that $ l_1 \geqslant \dots \geqslant l_p $ and
$l_1 + \dots + l_p \leqslant n $.

Answers to this problem can be found by counting the strings, and this
is implemented in \textsc{C} as \texttt{findec.c} (Appendix
\ref{app:CCode}). Output from this program (for $ n =1, \dots, 100 $)
is presented in Appendix \ref{app:ECData}. Not surprisingly, $ e $
grows rapidly with $ n $.  (This is just for curiosity purposes; will
only apply the case $n = 3 $.)

For the case $ n = 3 $, there are $ 7 $ types of classes of points,
listed in Table \ref{tab:generators}. There could be a generator at the
origin, so that's one type, called type $ [ 0 ] $.  Generators on a
coordinate axis are of a second type, called type $ [ 1 ] $, in which
there are $ 6 $ members of each equivalence class. Generators of the
form $ ( \beta, \beta, 0 ) $, in a class of size $ 12 $, are of type
$[ 1, 1 ] $, etc.%
\footnote{
  Each of these generator types can be thought of in terms of
  geometrical arrangements of points on the surface of a unit sphere,
  projected from the vertices, edges and faces of the unit octahedron
  \cite{Lebedev:76}.
}
More generally, the generator in (\ref{eq:generator}) is of type
$[l_1,l_2, \dots, l_p ] $, and the complete set of types of generators
in 3D is included in Table \ref{tab:generators}.  This notation (from
\cite{Keast:79}), is simplified for the case $ n = 3 $ in
\cite{MantelRabinowitz:77}, where there are $ K_i $ generators of each
type, for $ i = 0, \dots, 6 $, and we shall use the latter notation.
We shall refer to our rule as having \emph{structure}
${\left\lbrace K_i \right\rbrace}_{i=0}^6 $ (an ordered set), usually
just written $\left\lbrace K_i \right\rbrace $.

\begin{table}[htbp]
  \centering
  \begin{tabular}{||c|c|c|ccl|c||}
    \hline\hline
    Class & Class & Number of Generators & & & & Class \\
    Number & Type &
    (\protect\cite{Keast:79} and \protect\cite{MantelRabinowitz:77}) &
    \multicolumn{3}{c|}{Names of Generators and Weights} & Size \\
    \hline\hline
    $ 0 $ & $ [ 0       ] $ & $ K [ 0 ] = K_0 $ &
    $ ~ ( 0, 0, 0 ) $, & $ o   $ & \mbox{if~} $ K_0 = 1 $
    & $ ~1 $ \\
    $ 1 $ & $ [ 1       ] $ & $ K [ 1 ] = K_1 $ &
    $ ( \alpha_i,   0,          0          ) $, &
    $ a_i $ & $ i = 1, \dots, K_1    $ & $ ~6 $ \\
    $ 2 $ & $ [ 2       ] $ & $ K [ 2 ] = K_2 $ &
    $ ( \beta_i,    \beta_i,    0          ) $, &
    $ b_i $ & $ i = 1, \dots, K_2    $ & $ 12 $ \\
    $ 3 $ & $ [ 1, 1    ] $ & $ K [ 1, 1 ] = K_3 $ &
    $ ( \gamma_i,   \delta_i,   0          ) $, &
    $ c_i $ & $ i = 1, \dots, K_3    $ & $ 24 $ \\
    $ 4 $ & $ [ 3       ] $ & $ K [ 3 ] = K_4 $ &
    $ ( \epsilon_i, \epsilon_i, \epsilon_i ) $, &
    $ d_i $ & $ i = 1, \dots, K_4    $ & $ ~8 $ \\
    $ 5 $ & $ [ 2, 1    ] $ & $ K [ 2, 1 ] = K_5 $ &
    $ ( \zeta_i,    \zeta_i,    \eta_i     ) $, &
    $ e_i $ & $ i = 1, \dots, K_5    $ & $ 24 $ \\
    $ 6 $ & $ [ 1, 1, 1 ] $ & $ K [ 1, 1, 1 ] = K_6 $ &
    $ ( \theta_i,   \mu_i,      \lambda_i  ) $, &
    $ f_i $ & $ i = 1, \dots, K_6    $ & $ 48 $ \\
    \hline\hline
  \end{tabular}
  \caption
  {
    Nomenclature of generators and weights for Gau{\ss} cubature in 3D.
  }
  \label{tab:generators}
\end{table}

~\\

Given a structure $ \left\lbrace K_i \right\rbrace $, there will be a
total of $ N $ points in our rule, given by:
\begin{equation}
  N
  =
  K_0 + 6 K_1 + 12 K_2 + 24 K_3 + 8 K_4 + 24 K_5 + 48 K_6.
  \label{eq:CostFunction}
\end{equation}
To construct Gau{\ss} rules, we proceed \emph{without} presupposing the
structure, instead attempting to find a structure such that the
conditions for Gau{\ss} cubature are satisfied, and $ N $ is
minimised.  Taking $ m $ (such that the desired degree is $ 2m + 1 $),
and a rule structure $ \left\lbrace K_i \right\rbrace $, we write down
a system of (moment) equations involving an appropriate set of
generators, based on the requirement that a Gau{\ss} rule exactly
integrates all $n $-variable polynomials of each degree up to $ 2m $.
It is in fact a sufficient requirement that we integrate exactly
($n$-variable) monomials of these degrees.  The system of moment
equations is a system of nonlinear algebraic equations. (This is also
true in 1D, although clever artifice allows us to reduce its solution
to that of a linear system.)

For FS $\mathcal{R}_3 $, the variables are listed in Table
\ref{tab:generators}, and the system (\ref{eq:eval}) expands to that
presented in Figure \ref{fig:Star}.  We will call this system $(*)$, in
accordance with \cite[pp~410-411]{MantelRabinowitz:77}, where it first
appears explicitly.  System $(*)$ splits naturally into three
subsystems, defined by the number of non-zero indices $ k_i $ in the
monomial
$
  x_1^{2 k_1}
  x_2^{2 k_2}
  x_3^{2 k_3}
$
that we wish to integrate exactly.  Recall that the
$ \left\lbrace K_i\right\rbrace $ are non-negative integers, and the
other variables are real. Also, although
$K_0 \in \left\lbrace 0, 1 \right\rbrace $, we include it in a sum for
consistency, and apply the convention that $ \sum_{i=1}^{K_0} o = 0 $
if $ K_0 = 0 $.  If $ m < 3 $, subsystem III is ignored, and if $m<2$,
subsystem II is also ignored.

\begin{figure}[htbp]
  \setbox4=\vbox{
  \begin{eqnarray*}
    \mbox{Subsystem I:}
    \\
    I [ 1 ]
    & = &
       \sum_{i=1}^{K_0} o
    +
     6 \sum_{i=1}^{K_1} a_i
    +
    12 \sum_{i=1}^{K_2} b_i
    +
    24 \sum_{i=1}^{K_3} c_i
    +
     8 \sum_{i=1}^{K_4} d_i
    +
    24 \sum_{i=1}^{K_5} e_i
    +
    48 \sum_{i=1}^{K_6} f_i
    \\
    I [ x^{2j} ]
    & = &
     2 \sum_{i=1}^{K_1} a_i \alpha_i^{2j}
    +
     8 \sum_{i=1}^{K_2} b_i \beta_i^{2j}
    +
     8 \sum_{i=1}^{K_3} c_i ( \gamma_i^{2j} + \delta_i^{2j} )
    +
     8 \sum_{i=1}^{K_4} d_i \epsilon_i^{2j}
    \\
    & & +
     8 \sum_{i=1}^{K_5} e_i ( 2 \zeta_i^{2j} + \eta_i^{2j} )
    +
    16 \sum_{i=1}^{K_6}
      f_i
      ( \theta_i^{2j} + \mu_i^{2j} + \lambda_i^{2j} )
    \\
    & &
    \qquad \qquad \qquad \qquad \qquad \qquad
    j = 1, \dots, m
    \\
    \mbox{Subsystem II:}
    \\
    I [ x^{2j} y^{2k} ]
    & = &
     4 \sum_{i=1}^{K_2} b_i \beta_i^{2j + 2k}
    +
     4 \sum_{i=1}^{K_3}
      c_i
      (
        \gamma_i^{2j} \delta_i^{2k}
        +
        \gamma_i^{2k} \delta_i^{2j}
      )
    +
     8 \sum_{i=1}^{K_4} d_i \epsilon_i^{2j + 2k}
    \\
    & & +
     8 \sum_{i=1}^{K_5}
      e_i
      (
        \zeta_i^{2j + 2k}
        +
        \zeta_i^{2j} \eta_i^{2k}
        +
        \zeta_i^{2k} \eta_i^{2j}
      )
    \\
    & & +
     8 \sum_{i=1}^{K_6}
      f_i
      (
        \theta_i^{2j} \mu_i^{2k}
        +
        \theta_i^{2k} \mu_i^{2j}
        +
        \theta_i^{2j} \lambda_i^{2k}
        +
        \theta_i^{2k} \lambda_i^{2j}
        +
        \mu_i^{2j} \lambda_i^{2k}
        +
        \mu_i^{2k} \lambda_i^{2j}
      )
    \\
    & &
    \qquad \qquad \qquad \qquad \qquad \qquad
    1 \leqslant j \leqslant k
    \qquad
    j + k = 2, \dots, m
    \\
    \mbox{Subsystem III:}
    \\
    I [ x^{2j} y^{2k} z^{2l} ]
    & = &
     8 \sum_{i=1}^{K_4} d_i \epsilon_i^{2j + 2k + 2l}
    +
     8 \sum_{i=1}^{K_5}
      e_i
      (
        \zeta_i^{2j + 2k} \eta_i^{2l}
        +
        \zeta_i^{2j + 2l} \eta_i^{2k}
        +
        \zeta_i^{2k + 2l} \eta_i^{2j}
      )
    \\
    & &
    \! \! \! \! \! \! \! \! \! \! \! \! \! \! \! \!
    \! \! \! \! \! \! \! \! \! \! \! \! \! \! \! \!
    +
     8 \sum_{i=1}^{K_6}
      f_i
      (
        \theta_i^{2j} \mu_i^{2k} \lambda_i^{2l}
        +
        \theta_i^{2j} \mu_i^{2l} \lambda_i^{2k}
        +
        \theta_i^{2k} \mu_i^{2j} \lambda_i^{2l}
        +
        \theta_i^{2k} \mu_i^{2l} \lambda_i^{2j}
        +
        \theta_i^{2l} \mu_i^{2j} \lambda_i^{2k}
        +
        \theta_i^{2l} \mu_i^{2k} \lambda_i^{2j}
      )
    \\
    & &
    \qquad \qquad \qquad \qquad \qquad \qquad
    1 \leqslant j \leqslant k \leqslant l
    \qquad
    j + k + l = 3, \dots, m
  \end{eqnarray*}
  }
  \newdimen\xxx \xxx=\wd4 \advance\xxx by 60pt \wd4=\xxx
  \boxit{\boxit{\box4}}
  \caption{
    The system $(*)$ of moment equations used to determine
    Gau{\ss} cubature rules for regions
    $\mathcal{R}_3 \subseteq \protect\mathbb{R}^{3} $ (after
    \protect\cite[pp~410-411]{MantelRabinowitz:77}).
  }
  \label{fig:Star}
\end{figure}

Row-wise examination of the variables in Table \ref{tab:generators} shows
that, given $ \left\lbrace K_i \right\rbrace $, there are a total of:
\begin{eqnarray*}
  v
  =
  K_0 + 2 K_1 + 2 K_2 + 3 K_3 + 2 K_4 + 3 K_5 + 4 K_6
\end{eqnarray*}
variables in $(*)$. To determine the number of equations, let $ r $
be a positive integer and, for $ \nu = 1, \dots, n $, let
$p_\nu ( r ) $ be the number of solutions in positive integers
$k_i $ of:
\begin{eqnarray*}
  k_0 + k_1 + \dots + k_\nu
  =
  r
  \qquad \qquad
  1
  \leqslant
  k_0
  \leqslant
  k_1
  \leqslant
  \dots
  \leqslant
  k_\nu.
\end{eqnarray*}
For $ \nu = 0 $, we say there's one solution if $ r = 1 $, and none
otherwise.  $ p_{\nu} ( r ) $ is the dimension of the space spanned
by:
\begin{eqnarray*}
  \left\lbrace
    x_1^{2 k_1}
    x_2^{2 k_2}
    \cdots
    x_{\nu}^{2 k_{\nu}}
    \; | \;
    k_1
    = k_2
    \geqslant k_3
    \geqslant \cdots
    \geqslant k_{\nu}
    \geqslant 0
    \mbox{~and~}
    k_1 + k_2 + \dots + k_{\nu} = r
  \right\rbrace.
\end{eqnarray*}
This can be efficiently computed by:
\begin{eqnarray*}
  p_{\nu} ( r )
  & = &
  \left\lbrace
    \begin{array}{ll}
      \left\lbrace
        \begin{array}{ll}
          1 & r = 0 \\
          0 & \mbox{else}
        \end{array}
      \right.
                                 & \nu = 0 \\
      0                          & \nu = 1 \\
      r + 1 ( \mbox{mod~} 2 )  & \nu = 2 \\
      p_{\nu-1} ( r - 1 )
      +
      p_\nu ( r - \nu )        & \mbox{else}.
    \end{array}
  \right.
\end{eqnarray*}
For each $ r = 1, \dots, m $, there are $ p_0 ( r ) + p_1 ( r ) $
equations in subsystem I, $ p_2 ( r ) $ in II and $ p_3 ( r ) $ in
III.  Hence, for a rule of degree $ 2m + 1 $, the total number of
equations involved is $ \sum_{r=1}^m \sum_{\nu=0}^3 p_{\nu} ( r ) $.
For $ m = 1, \dots, 20 $, these numbers are listed in Table
\ref{tab:NumEqns}.

\begin{table}
  \centering
  \begin{tabular}{||c|*{4}{c}|c||}
    \hline\hline
    $ r $ & $ p_0 ( r ) $ & $ p_1 ( r ) $ & $ p_2 ( r ) $ &
    $ p_3 ( r ) $ &
    $
      \sum_{r=1}^m \sum_{\nu=0}^3
        p_{\nu} ( r )
    $
    \\
    \hline\hline
    ~1 &   1 &   1 &  ~0 &  ~0 & ~~2 \\
    ~2 &   0 &   1 &  ~1 &  ~0 & ~~4 \\
    ~3 &   0 &   1 &  ~1 &  ~1 & ~~7 \\
    ~4 &   0 &   1 &  ~2 &  ~1 & ~11 \\
    ~5 &   0 &   1 &  ~2 &  ~2 & ~16 \\
    ~6 &   0 &   1 &  ~3 &  ~3 & ~23 \\
    ~7 &   0 &   1 &  ~3 &  ~4 & ~31 \\
    ~8 &   0 &   1 &  ~4 &  ~5 & ~41 \\
    ~9 &   0 &   1 &  ~4 &  ~7 & ~53 \\
    10 &   0 &   1 &  ~5 &  ~8 & ~67 \\
    11 &   0 &   1 &  ~5 &  10 & ~83 \\
    12 &   0 &   1 &  ~6 &  12 & 102 \\
    13 &   0 &   1 &  ~6 &  14 & 123 \\
    14 &   0 &   1 &  ~7 &  16 & 147 \\
    15 &   0 &   1 &  ~7 &  19 & 174 \\
    16 &   0 &   1 &  ~8 &  21 & 204 \\
    17 &   0 &   1 &  ~8 &  24 & 237 \\
    18 &   0 &   1 &  ~9 &  27 & 274 \\
    19 &   0 &   1 &  ~9 &  30 & 314 \\
    20 &   0 &   1 &  10 &  33 & 358 \\
    \hline\hline
  \end{tabular}
  \caption
  {
    The total number of equations involved in system $(*)$,
    for $ m = 1, \dots, 20 $.
  }
  \label{tab:NumEqns}
\end{table}

\pagebreak

In general, $(*)$ may not have a unique solution, and indeed, may
have no solutions.  The best analytic technique currently available for
the solution of systems of multivariable polynomial equations is the use
of Groebner bases, see for example \cite{CoxLittleOShea:92}. (The
computer algebra package \textsc{Mathematica} uses this method.) Symbolic
computations are expensive, however, and will fail to yield answers for
polynomial equations of degree greater than $ 4 $.  We are quickly led
to numerical techniques!  Fortunately, experience with $(*)$
shows that solutions commonly \emph{do} exist; we progress with this as a
hope.

For a given degree $ 2m + 1 $, there are many different structures
$\left\lbrace K_i \right\rbrace $ that will lead to a system $(*)$.  To deduce a
possible structure, we would like to choose $ \left\lbrace K_i \right\rbrace $ to minimise
$N $, such that $(*)$ is consistent. This is an optimisation
problem, with constraints that will ensure the loosest possible
consistency of $(*)$. It turns out
\cite
{%
  KeastLyness:79,%
  MantelRabinowitz:77%
}
that the optimisation problem has \emph{linear} constraints. Clearly the
cost function (\ref{eq:CostFunction}) is linear, and the solutions
$\left\lbrace K_i \right\rbrace $ are integers.  The optimisation problem is thus a
linear integer programming problem, with linear cost function.  This is
routine to solve, and many examples of solutions are provided by
Rabinowitz et al. in \cite{RabinowitzRichter:69} (for 2D) and
\cite
{%
  MantelRabinowitz:74,%
  MantelRabinowitz:77%
} (for 3D).

\vfill

\pagebreak

We will call the constraints on $(*)$ \emph{consistency conditions};
their number will grow with dimension $ n $.  Establishment of the
consistency conditions is in general a formidable problem, and can only
be done manually for two or three dimensions.%
\footnote{
  Although a framework for the automatic construction of consistency
  conditions in any number of dimensions is presented in
  \cite{KeastLyness:79}, it seems to be essential to use a machine to
  do this.
}
For 3D, it initially appears that there could be as many as
$2^7 - 1 = 127 $ constraints, and in general, where $\mathcal{R}_n $ has
$e+1 $ types of FS equivalence classes, there might be $ 2^{e+1} - 1 $
constraints.  Fortunately, symmetries in $(*)$ reduce this number
to a manageable level.  For $ n = 3 $, this reduction can be done
manually,%
\footnote{
  See \cite[pp~394-398]{MantelRabinowitz:77}, and extra details in
  \cite{MantelRabinowitz:74}.
}
and leads to the system of $ 13 $ constraints presented in Figure \ref{fig:CC}.
The data presented in Table \ref{tab:NumEqns} can be used to convert these
constraints to a linear vector inequality.

\begin{figure}[htbp]
  \setbox4=\vbox{
  \begin{eqnarray*}
    3 K_5 + 4 K_6
    & \geqslant &
    \sum_{r=3}^m
      ( p_3 ( r ) - 1 )
    \\
    2 K_4 + 3 K_5 + 4 K_6
    & \geqslant &
    \sum_{r=3}^m
      p_3 ( r )
    \\
    3 K_3 + 3 K_5 + 4 K_6
    & \geqslant &
    \sum_{r=3}^m
      ( p_2 ( r ) + p_3 ( r ) - 2 )
    \\
    3 K_3 + 2 K_4 + 3 K_5 + 4 K_6
    & \geqslant &
    \sum_{r=2}^m
      ( p_2 ( r ) + p_3 ( r ) - 1 )
    \\
    2 K_2 + 3 K_3 + 3 K_5 + 4 K_6
    & \geqslant &
    \sum_{r=2}^m
      ( p_2 ( r ) + p_3 ( r ) - 1 )
    \\
    2 K_1 + 3 K_3 + 3 K_5 + 4 K_6
    & \geqslant &
    \sum_{r=2}^m
      ( p_2 ( r ) + p_3 ( r ) - 1 )
    \\
    2 K_2 + 3 K_3 + 2 K_4 + 3 K_5 + 4 K_6
    & \geqslant &
    \sum_{r=2}^m
      ( p_2 ( r ) + p_3 ( r ) )
    \\
    2 K_1 + 3 K_3 + 2 K_4 + 3 K_5 + 4 K_6
    & \geqslant &
    \sum_{r=2}^m
      ( p_2 ( r ) + p_3 ( r ) )
    \\
    2 K_1 + 2 K_2 + 3 K_3 + 3 K_5 + 4 K_6
    & \geqslant &
    \sum_{r=2}^m
      ( p_2 ( r ) + p_3 ( r ) )
    \\
    4 K_6
    & \geqslant &
    \sum_{r=9}^m
      ( p_3 ( r ) - ( r - 3 ) )
    \\
    3 K_3 + 4 K_6
    & \geqslant &
    \sum_{r=6}^m
      ( p_2 ( r ) + p_3 ( r ) - ( r - 1 ) )
    \\
    K_0 + 2 K_1 + 2 K_2 + 3 K_3 + 2 K_4 + 3 K_5 + 4 K_6
    & \geqslant &
    1
    +
    \sum_{r=1}^m
      ( p_1 ( r ) + p_2 ( r ) + p_3 ( r ) )
    \\
    K_0
    & \leqslant &
    1
  \end{eqnarray*}
  }
  \boxit{\boxit{\hspace{-8mm}\box4}}
  \caption
  {
    The constraints for system $(*)$, expressed in
    terms of the $ p_{\nu} ( r ) $ (after
    \protect\cite[p~398]{MantelRabinowitz:77}, in
    which it is called system $ \bar{C} $).
    We apply the convention that $ \sum_{i=p}^q a_i = 0 $
    if $ q < p $.
  }
  \label{fig:CC}
\end{figure}

Once the integer programming problem has been solved for a candidate
structure $ \left\lbrace K_i \right\rbrace $, we write the specific
version of $(*)$ for this structure, and attempt to solve it for
generators and weights.  The moments on the LHS of $(*)$ are
available analytically for a wide range of FS $\mathcal{R}_3 $ (in
particular $ U_3 $).  There are a number of subtleties in this
process:

\begin{enumerate}
\item
  Whilst the IPP will always have a solution, this is not guaranteed to
  be unique. We will in general have a (small) number of possible rule
  structures to choose from, and we order them lexically.
\item
  Unfortunately, $(*)$ is \emph{not} guaranteed to be consistent
  for any particular rule structure deduced from the conditions
  (remember that they were chosen as the loosest possible), so we may
  have to try a number of possible structures before we succeed.
\item
  Assuming that there is a solution for a particular rule structure,
  this is not guaranteed to be unique; typically there will be either a
  (small) finite, or an infinite number of solutions to $(*)$.
\item
  Rules deduced from solution of $(*)$ may not be good. They may
  still be acceptable, in cases where the points are not in
  $\mathcal{R}_n $, but they are `close' to it; or where some weights are
  negative, but small in magnitude. For purity, we reject such
  solutions.
\item
  If there are no rules for minimal $ N $, or those that do exist are
  not good, then we search for (FSGM) rules of a larger $ N $. To do
  this, we restart the IPP with an added constraint that $ N $ be
  larger than the rejected solution. The new IPP will be solvable, and
  we continue with the solution of $(*)$, iterating this
  procedure until we arrive at a good rule.  This process is guaranteed
  to terminate, the worst possible case being that we actually
  construct a product rule (which always exist).
\item
  The numerical solution of $(*)$ may be computationally
  intractable. It may be difficult to tell whether an algorithm is not
  converging because a solution does not exist, or because the system
  is so severely nonlinear that the software cannot handle it.
\end{enumerate}

Despite these problems, solutions to $(*)$ have been computed for
several $\mathcal{R}_n $, and various weighting functions. Primary
results for $\mathcal{R}_2 $ and $\mathcal{R}_3 $ are contained in
\cite
{%
  MantelRabinowitz:74,%
  MantelRabinowitz:77,%
  RabinowitzRichter:69%
}.
Some results for $ U_3 $ (see \S\ref{sec:U3Cubature}), are presented in
Appendix \ref{app:GaussRulesU3}.

\vfill

\pagebreak


\section{Gau{\ss} Cubature for $ U_3 $}
\label{sec:U3Cubature}

We specialise the theory presented in \S\ref{sec:GCforFS} to the case
where $\mathcal{R}_n $ is $ U_3 $, and the weighting function is unity.
As the `volume' of $ U_3 $ is finite, integrals over it of any bounded
function will be defined.  We describe the spherical harmonic
polynomials, which our cubature rule should integrate exactly, and
comment that this is equivalent to our previous requirement that our
rules integrate polynomials exactly. The construction of Gau{\ss}
cubature rules follows, with simplifying assumptions that are obtained
by considering the region of integration. The material presented is
related to that of Keast
\cite
{%
  KeastDiaz:83,%
  Keast:87%
}.

\subsection{Spherical Harmonics}

Spherical harmonic polynomials are the natural generalisation of
orthogonal polynomials to the surface of the sphere. They satisfy
Laplace's equation~\cite{Hobson:31}, appearing in particular when
separation of variables is used to solve the interior Dirichlet
problem. They are orthogonal, and hence are a useful basis for
approximating functions on $ U_3 $. From
\cite{PressFlanneryTeukolskyVetterling:88}, we take the following
definitions and properties.  Any point $ \mathbf{x} \in U_3 $ can be uniquely
characterised by the ordered pair $ ( \theta, \phi ) $ of coordinates
of longitude $ \theta $ and colatitude $ \phi $.  The
\emph{spherical harmonic polynomials} $ Y_{l m} ( \theta, \phi ) $,
where $ m \leqslant | l | $, and $ l, m \in \mathbb{Z} $ are defined in terms
of the \emph{associated Legendre polynomials} $ P_l^m $ as follows:
\begin{eqnarray*}
  Y_{l m} ( \theta, \phi )
  =
  \sqrt
  {
    {\displaystyle\frac
    {
      ( 2l + 1 ) ( l - m )!
    }{
      4 \pi ( l + m )!
    }}
  }
  \,
  P_l^m ( \cos ( \theta ) )
  e^{i m \theta}.
\end{eqnarray*}
They are orthogonal functions, normalised such that the integral over
$U_3 $ of a product $ Y_{l m} Y_{l' m'} $ is unity only if $ l = l' $
and $ m = m' $. Using $ * $ to denote complex conjugation:
\begin{eqnarray*}
  \int_0^{2 \pi}
    \!
    \int_{-1}^{1}
      Y_{l' m'}^{*} ( \theta, \phi )
      Y_{l m} ( \theta, \phi )
      \,
    d ( \cos ( \theta ) )
    \,
  d\phi
  =
  {\delta}_{l' l}
  {\delta}_{m' m}.
\end{eqnarray*}
Using the relation
\begin{eqnarray*}
  Y_{l,-m} ( \theta, \phi )
  =
  {( - )}^m
  Y_{l,-m}^{*} ( \theta, \phi ),
\end{eqnarray*}
we can relate any spherical harmonic to an associated Legendre
polynomial $ P_l^m $ with $ m \geqslant 0 $.  Table \ref{tab:SH} lists
some of the simplest spherical harmonics. As $ l $ and $ m $ increase,
the degree of the trigonometric functions increases. Note that for
$m=0$, they are purely real.

\renewcommand{\arraystretch}{2.1}
\begin{table}[htbp]
  \centering
  \begin{tabular}{||c|c|c||}
    \hline\hline
    $ l \quad m $  & $ P_l^m ( x ) $ &
    $ Y_{l m} ( \theta, \phi ) $ \\
    \hline\hline
    $ 0 \quad 0 $  & $ 1 $ &
    $ \sqrt{\frac{1}{4 \pi}} $ \\
    \hline
    $ 1 \quad 0 $  & $ x $ &
    $ \sqrt{\frac{3}{4 \pi}} \cos ( \theta ) $ \\
    $ 1 \quad 1 $  & $ - \sqrt{1 - x^2} $ &
    $ - \sqrt{\frac{3}{8 \pi}}
      \sin ( \theta ) e^{i \phi} $ \\
    \hline
    $ 2 \quad 0 $  & $ \frac{1}{2} ( 3 x^2 - 1 ) $ &
    $
      \sqrt{\frac{5}{4 \pi}}
      (
        \frac{3}{2} \cos^2 ( \theta )
        -
        \frac{1}{2}
      )
    $ \\
    $ 2 \quad 1 $  & $ - 3 x \sqrt{1 - x^2} $ &
    $
      - \sqrt{\frac{15}{8 \pi}}
      \sin ( \theta )
      \cos ( \theta )
      e^{i \phi}
    $ \\
    $ 2 \quad 2 $  & $ 3 ( 1 - x^2 ) $ &
    $
      \frac{1}{4}
      \sqrt{\frac{15}{2 \pi}}
      \sin^2 ( \theta )
      e^{2 i \phi}
    $ \\
    \hline\hline
  \end{tabular}
  \caption
  {
    The first few spherical harmonics and the associated Legendre
    polynomials
    (after \protect\cite[p~195]{PressFlanneryTeukolskyVetterling:88}).
  }
  \label{tab:SH}
\end{table}
\renewcommand{\arraystretch}{1}

\vfill

\pagebreak

Any function $ F : U_3 \to \mathbb{R} $ can be expanded in terms of
spherical harmonics:
\begin{eqnarray*}
  F ( \theta, \phi )
  =
  \sum_{l=0}^{\infty} \sum_{m=0}^{\infty}
    {\alpha}_{l m}
    Y_{l m} ( \theta, \phi ).
\end{eqnarray*}
Gau{\ss} quadrature rules in 1D exactly integrate monomials of as high
a degree as possible using a fixed number of quadrature points. In the
case of the cubature of integrals defined on $ U_3 $, the natural
generalisation is that the rules must exactly integrate as many
spherical harmonics as possible.%
\footnote{
  The 1D analogue of $ U_3 $ is the unit circle $ U_2 $, for which the
  orthogonal polynomials are $ e^{i n \phi} $ and the points of
  Gau{\ss} rules are equally spaced.
}
As in the 1D case, the cubature rule may not be optimal for any
particular integrand, but should work well for smooth functions that
are well approximated by (sums of) spherical harmonics.  It turns out
that exact integration of the spherical harmonics over $ U_3 $ is
equivalent to applying the work in \S\ref{sec:GCforFS}, with
restrictions on the placement of points. We do not have to consider
$U_3$ as a special case, and continue on from \S\ref{sec:GCin3D}.


\subsection{Application of the Gau{\ss} Cubature}

For $ U_3 $, the points of a good cubature rule (Table
\ref{tab:generators}) must lie on $ U_3 $, and so have some further
constraints added to them.  In particular $ K_0 = 0 $ and
$K_1,K_2,K_4\in \left\lbrace 0, 1 \right\rbrace $. Also, if $K_1=1$
then $\alpha_1 = 1 $, if $ K_2 = 1 $ then $ \beta_1 = 1/\sqrt{2} $, and
if $K_4 = 1 $ then $ \epsilon_1 = 1/\sqrt{3} $.  Furthermore:
\begin{eqnarray*}
  \begin{array}{ll}
    \gamma_i^2 + \delta_i^2 = 1
    &
    i = 1, \dots, K_3
    \\
    2 \zeta_i^2 + \eta_i^2 = 1
    &
    i = 1, \dots, K_5
    \\
    \theta_i^2 + \mu_i^2 + \lambda_i^2 = 1
    &
    i = 1, \dots, K_6.
  \end{array}
\end{eqnarray*}
In general for FS regions $\mathcal{R}_3 $, there are a total of $ 13 $
consistency conditions (Figure \ref{fig:CC}), but for $ U_3 $ these
simplify drastically. There are only $ 4 $, and these listed in Figure
\ref{fig:U3CC}.  For $ m = 1, \dots, 20 $, the appropriate right hand
sides of the constraint equations in Figure \ref{fig:U3CC} are
presented in Table \ref{tab:RHSU3}.

\begin{figure}[htbp]
  \setbox4=\vbox{
  \begin{eqnarray*}
    K_1 + K_2 + 2 K_3 + K_4 + 2 K_5 + 3 K_6
    & \geqslant &
    1 +
    \sum_{r=2}^m
      [ p_2 ( r ) + p_3 ( r ) ]
    \\
    K_4 + 2 K_5 + 3 K_6
    & \geqslant &
    1 +
    \sum_{r=3}^m
      p_3 ( r )
    \\
    2 K_3 + 3 K_6
    & \geqslant &
    1 +
    \sum_{r=6}^m
      [ p_2 ( r ) + p_3 ( r )  - 1 ]
    \\
    3 K_6
    & \geqslant &
    \sum_{r=9}^m
      [ p_3 ( r )  - 1 ]
    \end{eqnarray*}
  }
  \boxit{\boxit{\box4}}
  \caption
  {
    Constraints for the integer programming problem, where
    $\mathcal{R}_n $ is $ U_3 $ (c.f. Figure \protect\ref{fig:CC}). In addition,
    $ K_0 = 0 $, and $ K_1, K_2, K_4 \leqslant 1 $.  Again, the
    convention $ \sum_{i=p}^q a_i = 0 $ if $ q < p $, is followed.
    Values for the right hand sides are tabulated for some choices of
    $ r $ in Table \protect\ref{tab:RHSU3}.
  }
  \label{fig:U3CC}
\end{figure}

\begin{table}[htbp]
  \small
  \centering
  \begin{tabular}{||c|*{20}{r}||}
    \hline\hline
    $ m $
    &  1 &  2 &  3 &  4 &  5 &  6 &  7 &  8 &  9 & 10 &
      11 & 12 & 13 & 14 & 15 & 16 & 17 & 18 & 19 & 20 \\
    \hline\hline
    &  1 &  2 &  3 &  4 &  5 &  7 &  8 & 10 & 12 & 14 &
      16 & 19 & 21 & 24 & 27 & 30 & 33 & 37 & 40 & 44 \\
    & 0 &  0 &  1 &  1 &  2 &  3 &  4 &  5 &  7 &  8 &
     10 & 12 & 14 & 16 & 19 & 21 & 24 & 27 & 30 & 33 \\
    & 0 &  0 &  0 &  0 &  0 &  1 &  1 &  2 &  3 &  4 &
      5 &  7 &  8 & 10 & 12 & 14 & 16 & 19 & 21 & 24 \\
    & 0 &  0 &  0 &  0 &  0 &  0 &  0 &  0 &  1 &  1 &
      2 &  3 &  4 &  5 &  7 &  8 & 10 & 12 & 14 & 16 \\
    \hline\hline
  \end{tabular}
  \normalsize
  \caption
  {
    The four elements of the columns of the right hand side in
    Figure \protect\ref{fig:U3CC}, for various $ m $, as generated by
    appropriate sums of the $ p_{\nu} ( r ) $.
  }
  \label{tab:RHSU3}
\end{table}

\pagebreak

Given $ m $, for which we wish to construct a rule of degree $2m + 1 $,
we firstly set up the integer programming problem to be solved for a
rule structure $ \left\lbrace K_i \right\rbrace $.  As in Table
\ref{tab:NumEqns}, the sum of the $ p_{\nu} ( r ) $ over
$ \nu = 0, \dots, 3 $, and $r = 1, \dots, m $ gives the total number of
equations in $(*)$ for $ U_3 $.

Lebedev \cite{Lebedev:76} published an important paper on Gau{\ss}
cubature for $ U_3 $ just prior to that of the more general one by
Mantel and Rabinowitz \cite{MantelRabinowitz:77}. Lebedev makes some
astute choices in the presupposition of rule structures, which
sometimes result in optimal choices. The algebra is simplified by
enforcing $ K_1 = K_2 = K_4 = 1 $, and sometimes also $ K_2 = 0 $; and
choosing $ K_3, K_4 $, and $ K_5 $ such that the number of unknowns in
$( * ) $ is equal to the number of equations.  The rules generated
were a great improvement on the previously completely unsystematised
collection of known rules (e.g. see Stroud \cite{Stroud:71}), but we
are interested in the more general case.


\begin{table}[htbp]
  \centering
  \begin{tabular}{||c|c||*{6}{c|}|}
    \hline\hline
    $ m $ & $ N $ &
    $ K_1 $ & $ K_2 $ & $ K_3 $ & $ K_4 $ & $ K_5 $ & $ K_6 $ \\
    \hline\hline
    $ ~1 $ & $ ~~6 $ & $ 1 $ & $ 0 $ & $ 0 $ & $ 0 $ & $ 0 $ & $ 0 $ \\
    $ ~2 $ & $ ~14 $ & $ 1 $ & $ 0 $ & $ 0 $ & $ 1 $ & $ 0 $ & $ 0 $ \\
    $ ~3 $ & $ ~26 $ & $ 1 $ & $ 1 $ & $ 0 $ & $ 1 $ & $ 0 $ & $ 0 $ \\
    $ ~4 $ & $ ~38 $ & $ 1 $ & $ 1 $ & $ 0 $ & $ 1 $ & $ 1 / 2 $ & $ 0 $ \\
    $ ~5 $ & $ ~50 $ & $ 1 $ & $ 1 $ & $ 0 $ & $ 1 $ & $ 1 $ & $ 0 $ \\
    $ ~6 $ & $ ~74 $ & $ 1 $ & $ 1 $ & $ 1 / 2 $ & $ 1 $ & $ 3/2 $ & $ 0 $ \\
    $ ~7 $ & $ ~86 $ & $ 1 $ & $ 1 $ & $ 1 / 2 $ & $ 1 $ & $ 2 $ & $ 0 $ \\
    $ ~8 $ & $ 110 $ & $ 1 $ & $ 1 $ & $ 1 $ & $ 1 $ & $ 5 / 2 $ & $ 0 $ \\
    $ ~9 $ & $ 138 $ & $ 1 $ & $ 1 $ & $ 1 $ & $ 1 $ & $ 3 $ & $ 1 / 3 $ \\
    $ 10 $ & $ 162 $ & $ 1 $ & $ 1 $ & $ 3 / 2 $ & $ 1 $ & $ 7/2 $ & $ 1/3 $ \\
    $ 11 $ & $ 190 $ & $ 1 $ & $ 1 $ & $ 3 / 2 $ & $ 1 $ & $ 4 $ & $ 2/3 $ \\
    $ 12 $ & $ 230 $ & $ 1 $ & $ 1 $ & $ 2 $ & $ 1 $ & $ 9 / 2 $ & $ 1 $ \\
    $ 13 $ & $ 258 $ & $ 1 $ & $ 1 $ & $ 2 $ & $ 1 $ & $ 5 $ & $ 4 / 3 $ \\
    $ 14 $ & $ 298 $ & $ 1 $ & $ 1 $ & $ 5 / 2 $ & $ 1 $ & $ 11/2 $ & $ 5/3 $ \\
    $ 15 $ & $ 342 $ & $ 1 $ & $ 1 $ & $ 5 / 2 $ & $ 1 $ & $ 6 $ & $ 7/3 $ \\
    $ 16 $ & $ 382 $ & $ 1 $ & $ 1 $ & $ 3 $ & $ 1 $ & $ 13 / 2 $ & $ 8/3 $ \\
    $ 17 $ & $ 426 $ & $ 1 $ & $ 1 $ & $ 3 $ & $ 1 $ & $ 7 $ & $ 10 / 3 $ \\
    $ 18 $ & $ 482 $ & $ 1 $ & $ 1 $ & $ 7 / 2 $ & $ 1 $ & $ 15/2 $ & $ 4 $ \\
    $ 19 $ & $ 526 $ & $ 1 $ & $ 1 $ & $ 7 / 2 $ & $ 1 $ & $ 8 $ & $ 14/3 $ \\
    $ 20 $ & $ 582 $ & $ 1 $ & $ 1 $ & $ 4 $ & $ 1 $ & $ 17 / 2 $ & $ 16/3 $ \\
    \hline\hline
  \end{tabular}
  \caption
  {
    Basic solutions to the IPP, with the integer constraint removed.
    Any integer solutions will have $ N $ at least equal to that
    displayed here.
  }
  \label{tab:IPPBasicSoln}
\end{table}

\pagebreak

\subsubsection{Solution of the Integer Programming Problem for $ U_3 $}

For our chosen $ m $, we set up and solve the integer programming
problem (including the constraints $ K_1, K_2, K_4 \leqslant 1 $, and
that $ K_0 $ is set to $ 0 $), for the minimisation of $ N $, the
number of points in the rule:
\begin{eqnarray*}
  N
  =
  6 K_1 + 12 K_2 + 24 K_3 + 8 K_4 + 24 K_5 + 48 K_6.
\end{eqnarray*}
We require \emph{all\/} the solutions to the IPP which have $ N $ at
least equal to the (integer) minimum, and less than some
(user-specified) small multiple of this minimum.  We firstly solve the
programming problem without integer constraints.  This yields a lower
bound for $N $, used when searching for \emph{all} the integer
solutions. This initial problem is solved using the
\textsc{Mathematica} program \texttt{IPPBasicSoln.M} (Appendix
\ref{app:MathematicaCode}).  Results (lower bounds for $ N $, and
associated `pseudostructures' $ \left\lbrace K_i \right\rbrace $), for
$ m = 1, \dots, 20 $ are presented in Table \ref{tab:IPPBasicSoln}.

To find \emph{all} the solutions of the IPP, we use a simple exhaustive
search. Whilst this is crude, in this case it is reasonably efficient.
For $ U_3 $, there are only six variables to search through, and
$K_1, K_2, K_4 \in \left\lbrace 0, 1 \right\rbrace $. We use the lower
bound on $ N $, and as an upper bound we use a small multiple (say
$1.5$) of the lower bound.  This second phase is implemented as a
\textsc{C} program \texttt{ipp.c} (Appendix \ref{app:CCode}). This
program takes the data in Table \ref{tab:IPPBasicSoln}, and
exhaustively searches for all integer solutions, subject to
(reasonable) bounds $K_3, K_5, K_6 \leqslant 20 $, and an upper bound
on $ N $ set empirically so that we only collect about the first $100$
solutions.  The output from \texttt{ipp.c} is (manually) sorted, so as
to order the solutions firstly with increasing $ N $, and then
lexically.  The solutions corresponding to the first five integer
minima, for $m = 1, \dots, 10 $ are presented in Table
\ref{tab:IPPTrueSoln}.

These minima agree with those published by Keast
\cite[p~155~and~pp~166-167]{Keast:87}. (This paper presents all
structures for the first $ 5 $ consecutive minima in $ N $, for
$m = 1, \dots, 9 $, corresponding to degrees $ 3, 5, \dots, 17 $.)
Keast claims to have obtained FSMG rules from the structures, but does
not actually describe them, but we do. (Keast's paper also correctly
identifies an error in the work of Lebedev \cite[p~15]{Lebedev:76}, but
fails to note that the Lebedev's paper is less general, so that results
cannot be compared directly.)

Following the notation of \cite[p~400]{MantelRabinowitz:77}, the rules
are assigned names. Say we have a rule of degree $ 2m + 1 $ for
$\mathcal{R}_n $ with weighting function $ h ( r ) $.  Let this rule be
found from the $ j $th instance (ordered lexically) of the $ i $th
consecutive minima of the IPP, and have structure
$\left\lbrace K_0, K_1, \dots, K_e \right\rbrace $, and total number of
points $ N $. We will name this rule:
\begin{eqnarray*}
 \mathcal{R}_n^{h ( r )} :
  ( 2m+1 ) \mbox{--} i.j ( K_0, K_1, \dots, K_e ) \mbox{--} N.
\end{eqnarray*}
Our rules for $ U_3 $, with $ h \equiv 1 $ can be labelled as:
\begin{eqnarray*}
  U_3 : ( 2m+1 ) \mbox{--} i.j ( K_1, \dots, K_6 ) \mbox{--} N.
\end{eqnarray*}
When $(*)$ has multiple solutions, the labelling is not unique,
but this is not a problem.  This nomenclature provides a basis for
comparison with (published) rules derived by other means.

\pagebreak

\tablehead
{
  \hline\hline\multicolumn{11}{||l||}%
  {\small\sl \dots continued from previous page}\\
  \hline
  $ m $ & $ i $ & $ j $ & $ N $ &
  $ K_1 $ & $ K_2 $ & $ K_3 $ & $ K_4 $ & $ K_5 $ & $ K_6 $ &
  $ v $ \\
  \hline\hline
}
\tabletail
{
  \hline\multicolumn{11}{||l||}%
  {\small\sl \dots continued on next page}\\
  \hline\hline
}
\tablelasttail{\hline\hline}
\tablefirsthead
{
  \hline\hline
  $ m $ & $ i $ & $ j $ & $ N $ &
  $ K_1 $ & $ K_2 $ & $ K_3 $ & $ K_4 $ & $ K_5 $ & $ K_6 $ &
  $ v $ \\
  \hline\hline
}
\bottomcaption
{
  All the solutions to the IPP, for $ m = 1, \dots, 10 $, corresponding
  to the first five minima. (Some of these are already available in
  Table \protect\ref{tab:IPPBasicSoln}.) $ v $ is the number of variables in
  $(*)$, for the particular structure.
}
{\centering
  \label{tab:IPPTrueSoln}
  \begin{supertabular}{||c|c|c||c||cccccc|c||}
    $  1 $&$ 1 $&$ 1 $&$   6 $&$ 1 $&$ 0 $&$  0 $&$ 0 $&$  0 $&$  0 $&$  2 $ \\
    $  1 $&$ 2 $&$ 1 $&$   8 $&$ 0 $&$ 0 $&$  0 $&$ 1 $&$  0 $&$  0 $&$  2 $ \\
    $  1 $&$ 3 $&$ 1 $&$  12 $&$ 0 $&$ 1 $&$  0 $&$ 0 $&$  0 $&$  0 $&$  2 $ \\
    $  1 $&$ 4 $&$ 1 $&$  14 $&$ 1 $&$ 0 $&$  0 $&$ 1 $&$  0 $&$  0 $&$  4 $ \\
    $  1 $&$ 5 $&$ 1 $&$  18 $&$ 1 $&$ 1 $&$  0 $&$ 0 $&$  0 $&$  0 $&$  4 $ \\
    \hline
    $  2 $&$ 1 $&$ 1 $&$  14 $&$ 1 $&$ 0 $&$  0 $&$ 1 $&$  0 $&$  0 $&$  4 $ \\
    $  2 $&$ 2 $&$ 1 $&$  18 $&$ 1 $&$ 1 $&$  0 $&$ 0 $&$  0 $&$  0 $&$  4 $ \\
    $  2 $&$ 3 $&$ 1 $&$  20 $&$ 0 $&$ 1 $&$  0 $&$ 1 $&$  0 $&$  0 $&$  4 $ \\
    $  2 $&$ 4 $&$ 1 $&$  24 $&$ 0 $&$ 0 $&$  0 $&$ 0 $&$  1 $&$  0 $&$  3 $ \\
    $  2 $&$ 4 $&$ 2 $&$  24 $&$ 0 $&$ 0 $&$  1 $&$ 0 $&$  0 $&$  0 $&$  3 $ \\
    $  2 $&$ 5 $&$ 1 $&$  26 $&$ 1 $&$ 1 $&$  0 $&$ 1 $&$  0 $&$  0 $&$  6 $ \\
    \hline
    $  3 $&$ 1 $&$ 1 $&$  26 $&$ 1 $&$ 1 $&$  0 $&$ 1 $&$  0 $&$  0 $&$  6 $ \\
    $  3 $&$ 2 $&$ 1 $&$  30 $&$ 1 $&$ 0 $&$  0 $&$ 0 $&$  1 $&$  0 $&$  5 $ \\
    $  3 $&$ 3 $&$ 1 $&$  32 $&$ 0 $&$ 0 $&$  0 $&$ 1 $&$  1 $&$  0 $&$  5 $ \\
    $  3 $&$ 3 $&$ 2 $&$  32 $&$ 0 $&$ 0 $&$  1 $&$ 1 $&$  0 $&$  0 $&$  5 $ \\
    $  3 $&$ 4 $&$ 1 $&$  36 $&$ 0 $&$ 1 $&$  0 $&$ 0 $&$  1 $&$  0 $&$  5 $ \\
    $  3 $&$ 5 $&$ 1 $&$  38 $&$ 1 $&$ 0 $&$  0 $&$ 1 $&$  1 $&$  0 $&$  7 $ \\
    $  3 $&$ 5 $&$ 2 $&$  38 $&$ 1 $&$ 0 $&$  1 $&$ 1 $&$  0 $&$  0 $&$  7 $ \\
    $  4 $&$ 1 $&$ 1 $&$  38 $&$ 1 $&$ 0 $&$  0 $&$ 1 $&$  1 $&$  0 $&$  7 $ \\
    $  4 $&$ 1 $&$ 2 $&$  38 $&$ 1 $&$ 0 $&$  1 $&$ 1 $&$  0 $&$  0 $&$  7 $ \\
    $  4 $&$ 2 $&$ 1 $&$  42 $&$ 1 $&$ 1 $&$  0 $&$ 0 $&$  1 $&$  0 $&$  7 $ \\
    $  4 $&$ 3 $&$ 1 $&$  44 $&$ 0 $&$ 1 $&$  0 $&$ 1 $&$  1 $&$  0 $&$  7 $ \\
    $  4 $&$ 3 $&$ 2 $&$  44 $&$ 0 $&$ 1 $&$  1 $&$ 1 $&$  0 $&$  0 $&$  7 $ \\
    $  4 $&$ 4 $&$ 1 $&$  48 $&$ 0 $&$ 0 $&$  0 $&$ 0 $&$  2 $&$  0 $&$  6 $ \\
    $  4 $&$ 4 $&$ 2 $&$  48 $&$ 0 $&$ 0 $&$  1 $&$ 0 $&$  1 $&$  0 $&$  6 $ \\
    $  4 $&$ 5 $&$ 1 $&$  50 $&$ 1 $&$ 1 $&$  0 $&$ 1 $&$  1 $&$  0 $&$  9 $ \\
    $  4 $&$ 5 $&$ 2 $&$  50 $&$ 1 $&$ 1 $&$  1 $&$ 1 $&$  0 $&$  0 $&$  9 $ \\
    \hline
    $  5 $&$ 1 $&$ 1 $&$  50 $&$ 1 $&$ 1 $&$  0 $&$ 1 $&$  1 $&$  0 $&$  9 $ \\
    $  5 $&$ 2 $&$ 1 $&$  54 $&$ 1 $&$ 0 $&$  0 $&$ 0 $&$  2 $&$  0 $&$  8 $ \\
    $  5 $&$ 2 $&$ 2 $&$  54 $&$ 1 $&$ 0 $&$  1 $&$ 0 $&$  1 $&$  0 $&$  8 $ \\
    $  5 $&$ 3 $&$ 1 $&$  56 $&$ 0 $&$ 0 $&$  0 $&$ 1 $&$  2 $&$  0 $&$  8 $ \\
    $  5 $&$ 3 $&$ 2 $&$  56 $&$ 0 $&$ 0 $&$  1 $&$ 1 $&$  1 $&$  0 $&$  8 $ \\
    $  5 $&$ 4 $&$ 1 $&$  60 $&$ 0 $&$ 1 $&$  0 $&$ 0 $&$  2 $&$  0 $&$  8 $ \\
    $  5 $&$ 4 $&$ 2 $&$  60 $&$ 0 $&$ 1 $&$  1 $&$ 0 $&$  1 $&$  0 $&$  8 $ \\
    $  5 $&$ 5 $&$ 1 $&$  62 $&$ 1 $&$ 0 $&$  0 $&$ 1 $&$  0 $&$  1 $&$  8 $ \\
    $  5 $&$ 5 $&$ 2 $&$  62 $&$ 1 $&$ 0 $&$  0 $&$ 1 $&$  2 $&$  0 $&$ 10 $ \\
    $  5 $&$ 5 $&$ 3 $&$  62 $&$ 1 $&$ 0 $&$  1 $&$ 1 $&$  1 $&$  0 $&$ 10 $ \\
    \hline
    $  6 $&$ 1 $&$ 1 $&$  74 $&$ 1 $&$ 1 $&$  1 $&$ 1 $&$  1 $&$  0 $&$ 12 $ \\
    $  6 $&$ 2 $&$ 1 $&$  78 $&$ 1 $&$ 0 $&$  1 $&$ 0 $&$  2 $&$  0 $&$ 11 $ \\
    $  6 $&$ 3 $&$ 1 $&$  80 $&$ 0 $&$ 0 $&$  1 $&$ 1 $&$  2 $&$  0 $&$ 11 $ \\
    $  6 $&$ 3 $&$ 2 $&$  80 $&$ 0 $&$ 0 $&$  2 $&$ 1 $&$  1 $&$  0 $&$ 11 $ \\
    $  6 $&$ 4 $&$ 1 $&$  84 $&$ 0 $&$ 1 $&$  1 $&$ 0 $&$  2 $&$  0 $&$ 11 $ \\
    $  6 $&$ 5 $&$ 1 $&$  86 $&$ 1 $&$ 0 $&$  0 $&$ 1 $&$  1 $&$  1 $&$ 11 $ \\
    $  6 $&$ 5 $&$ 2 $&$  86 $&$ 1 $&$ 0 $&$  1 $&$ 1 $&$  0 $&$  1 $&$ 11 $ \\
    $  6 $&$ 5 $&$ 3 $&$  86 $&$ 1 $&$ 0 $&$  1 $&$ 1 $&$  2 $&$  0 $&$ 13 $ \\
    $  6 $&$ 5 $&$ 4 $&$  86 $&$ 1 $&$ 0 $&$  2 $&$ 1 $&$  1 $&$  0 $&$ 13 $ \\
    \hline
    $  7 $&$ 1 $&$ 1 $&$  86 $&$ 1 $&$ 0 $&$  1 $&$ 1 $&$  2 $&$  0 $&$ 13 $ \\
    $  7 $&$ 2 $&$ 1 $&$  90 $&$ 1 $&$ 1 $&$  1 $&$ 0 $&$  2 $&$  0 $&$ 13 $ \\
    $  7 $&$ 3 $&$ 1 $&$  92 $&$ 0 $&$ 1 $&$  1 $&$ 1 $&$  2 $&$  0 $&$ 13 $ \\
    $  7 $&$ 4 $&$ 1 $&$  96 $&$ 0 $&$ 0 $&$  1 $&$ 0 $&$  3 $&$  0 $&$ 12 $ \\
    $  7 $&$ 4 $&$ 2 $&$  96 $&$ 0 $&$ 0 $&$  2 $&$ 0 $&$  2 $&$  0 $&$ 12 $ \\
    $  7 $&$ 5 $&$ 1 $&$  98 $&$ 1 $&$ 1 $&$  0 $&$ 1 $&$  1 $&$  1 $&$ 13 $ \\
    $  7 $&$ 5 $&$ 2 $&$  98 $&$ 1 $&$ 1 $&$  1 $&$ 1 $&$  0 $&$  1 $&$ 13 $ \\
    $  7 $&$ 5 $&$ 3 $&$  98 $&$ 1 $&$ 1 $&$  1 $&$ 1 $&$  2 $&$  0 $&$ 15 $ \\
    \hline
    $  8 $&$ 1 $&$ 1 $&$ 110 $&$ 1 $&$ 0 $&$  1 $&$ 1 $&$  3 $&$  0 $&$ 16 $ \\
    $  8 $&$ 1 $&$ 2 $&$ 110 $&$ 1 $&$ 0 $&$  2 $&$ 1 $&$  2 $&$  0 $&$ 16 $ \\
    $  8 $&$ 2 $&$ 1 $&$ 114 $&$ 1 $&$ 1 $&$  1 $&$ 0 $&$  3 $&$  0 $&$ 16 $ \\
    $  8 $&$ 3 $&$ 1 $&$ 116 $&$ 0 $&$ 1 $&$  1 $&$ 1 $&$  3 $&$  0 $&$ 16 $ \\
    $  8 $&$ 3 $&$ 2 $&$ 116 $&$ 0 $&$ 1 $&$  2 $&$ 1 $&$  2 $&$  0 $&$ 16 $ \\
    $  8 $&$ 4 $&$ 1 $&$ 120 $&$ 0 $&$ 0 $&$  1 $&$ 0 $&$  4 $&$  0 $&$ 15 $ \\
    $  8 $&$ 4 $&$ 2 $&$ 120 $&$ 0 $&$ 0 $&$  2 $&$ 0 $&$  3 $&$  0 $&$ 15 $ \\
    $  8 $&$ 5 $&$ 1 $&$ 122 $&$ 1 $&$ 1 $&$  0 $&$ 1 $&$  2 $&$  1 $&$ 16 $ \\
    $  8 $&$ 5 $&$ 2 $&$ 122 $&$ 1 $&$ 1 $&$  1 $&$ 1 $&$  1 $&$  1 $&$ 16 $ \\
    $  8 $&$ 5 $&$ 3 $&$ 122 $&$ 1 $&$ 1 $&$  1 $&$ 1 $&$  3 $&$  0 $&$ 18 $ \\
    $  8 $&$ 5 $&$ 4 $&$ 122 $&$ 1 $&$ 1 $&$  2 $&$ 1 $&$  2 $&$  0 $&$ 18 $ \\
    \hline
    $  9 $&$ 1 $&$ 1 $&$ 146 $&$ 1 $&$ 1 $&$  0 $&$ 1 $&$  3 $&$  1 $&$ 19 $ \\
    $  9 $&$ 1 $&$ 2 $&$ 146 $&$ 1 $&$ 1 $&$  1 $&$ 1 $&$  2 $&$  1 $&$ 19 $ \\
    $  9 $&$ 2 $&$ 1 $&$ 150 $&$ 1 $&$ 0 $&$  0 $&$ 0 $&$  4 $&$  1 $&$ 18 $ \\
    $  9 $&$ 2 $&$ 2 $&$ 150 $&$ 1 $&$ 0 $&$  1 $&$ 0 $&$  3 $&$  1 $&$ 18 $ \\
    $  9 $&$ 2 $&$ 3 $&$ 150 $&$ 1 $&$ 0 $&$  2 $&$ 0 $&$  2 $&$  1 $&$ 18 $ \\
    $  9 $&$ 3 $&$ 1 $&$ 152 $&$ 0 $&$ 0 $&$  0 $&$ 1 $&$  4 $&$  1 $&$ 18 $ \\
    $  9 $&$ 3 $&$ 2 $&$ 152 $&$ 0 $&$ 0 $&$  1 $&$ 1 $&$  3 $&$  1 $&$ 18 $ \\
    $  9 $&$ 3 $&$ 3 $&$ 152 $&$ 0 $&$ 0 $&$  2 $&$ 1 $&$  2 $&$  1 $&$ 18 $ \\
    $  9 $&$ 4 $&$ 1 $&$ 156 $&$ 0 $&$ 1 $&$  0 $&$ 0 $&$  4 $&$  1 $&$ 18 $ \\
    $  9 $&$ 4 $&$ 2 $&$ 156 $&$ 0 $&$ 1 $&$  1 $&$ 0 $&$  3 $&$  1 $&$ 18 $ \\
    $  9 $&$ 4 $&$ 3 $&$ 156 $&$ 0 $&$ 1 $&$  2 $&$ 0 $&$  2 $&$  1 $&$ 18 $ \\
    $  9 $&$ 5 $&$ 1 $&$ 158 $&$ 1 $&$ 0 $&$  0 $&$ 1 $&$  2 $&$  2 $&$ 18 $ \\
    $  9 $&$ 5 $&$ 2 $&$ 158 $&$ 1 $&$ 0 $&$  0 $&$ 1 $&$  4 $&$  1 $&$ 20 $ \\
    $  9 $&$ 5 $&$ 3 $&$ 158 $&$ 1 $&$ 0 $&$  1 $&$ 1 $&$  1 $&$  2 $&$ 18 $ \\
    $  9 $&$ 5 $&$ 4 $&$ 158 $&$ 1 $&$ 0 $&$  1 $&$ 1 $&$  3 $&$  1 $&$ 20 $ \\
    $  9 $&$ 5 $&$ 5 $&$ 158 $&$ 1 $&$ 0 $&$  2 $&$ 1 $&$  0 $&$  2 $&$ 18 $ \\
    $  9 $&$ 5 $&$ 6 $&$ 158 $&$ 1 $&$ 0 $&$  2 $&$ 1 $&$  2 $&$  1 $&$ 20 $ \\
    \hline
    $ 10 $&$ 1 $&$ 1 $&$ 170 $&$ 1 $&$ 1 $&$  1 $&$ 1 $&$  3 $&$  1 $&$ 22 $ \\
    $ 10 $&$ 1 $&$ 2 $&$ 170 $&$ 1 $&$ 1 $&$  2 $&$ 1 $&$  2 $&$  1 $&$ 22 $ \\
    $ 10 $&$ 2 $&$ 1 $&$ 174 $&$ 1 $&$ 0 $&$  1 $&$ 0 $&$  4 $&$  1 $&$ 21 $ \\
    $ 10 $&$ 2 $&$ 2 $&$ 174 $&$ 1 $&$ 0 $&$  2 $&$ 0 $&$  3 $&$  1 $&$ 21 $ \\
    $ 10 $&$ 3 $&$ 1 $&$ 176 $&$ 0 $&$ 0 $&$  1 $&$ 1 $&$  4 $&$  1 $&$ 21 $ \\
    $ 10 $&$ 3 $&$ 2 $&$ 176 $&$ 0 $&$ 0 $&$  2 $&$ 1 $&$  3 $&$  1 $&$ 21 $ \\
    $ 10 $&$ 3 $&$ 3 $&$ 176 $&$ 0 $&$ 0 $&$  3 $&$ 1 $&$  2 $&$  1 $&$ 21 $ \\
    $ 10 $&$ 4 $&$ 1 $&$ 180 $&$ 0 $&$ 1 $&$  1 $&$ 0 $&$  4 $&$  1 $&$ 21 $ \\
    $ 10 $&$ 4 $&$ 2 $&$ 180 $&$ 0 $&$ 1 $&$  2 $&$ 0 $&$  3 $&$  1 $&$ 21 $ \\
    $ 10 $&$ 5 $&$ 1 $&$ 182 $&$ 1 $&$ 0 $&$  0 $&$ 1 $&$  3 $&$  2 $&$ 21 $ \\
    $ 10 $&$ 5 $&$ 2 $&$ 182 $&$ 1 $&$ 0 $&$  1 $&$ 1 $&$  2 $&$  2 $&$ 21 $ \\
    $ 10 $&$ 5 $&$ 3 $&$ 182 $&$ 1 $&$ 0 $&$  1 $&$ 1 $&$  4 $&$  1 $&$ 23 $ \\
    $ 10 $&$ 5 $&$ 4 $&$ 182 $&$ 1 $&$ 0 $&$  2 $&$ 1 $&$  1 $&$  2 $&$ 21 $ \\
    $ 10 $&$ 5 $&$ 5 $&$ 182 $&$ 1 $&$ 0 $&$  2 $&$ 1 $&$  3 $&$  1 $&$ 23 $ \\
    $ 10 $&$ 5 $&$ 6 $&$ 182 $&$ 1 $&$ 0 $&$  3 $&$ 1 $&$  2 $&$  1 $&$ 23 $ \\
  \end{supertabular}
}

\pagebreak


\subsubsection{Moments for $ U_3 $}
\label{sec:U3Moments}

Solution of $(*)$ for $ U_3 $ requires knowledge of the moments
$I [ x^{2 j_1} y^{2 j_2} z^{2 j_3} ] $, for all combinations of
integers $ 0 \leqslant j_1 \leqslant j_2 \leqslant j_3 \leqslant m $,
such that $ 0 \leqslant j_1 + j_2 + j_3 \leqslant m $. (The moment is
zero if any of the exponents are odd, and it is invariant under
permutation of exponents.) Explicitly, we require analytical evaluation
of integrals of the form:
\begin{eqnarray*}
  I [ x^{2 j_1} y^{2 j_2} z^{2 j_3} ]
  =
  \int_{U_3}
    x^{2 j_1} y^{2 j_2} z^{2 j_3}
  \;
  dx \, dy \, dz.
\end{eqnarray*}
Clearly $ I [ x^0 y^0 z^0 ] = 4 \pi $.  Using coordinates $ \theta $,
the longitude, and $ \phi $, the latitude (not the \emph{co}-latitude),
the integral separates \cite[p~33]{Stroud:71}:
\begin{eqnarray*}
  \hspace{-9mm}
  I [ x^{2 j_1} y^{2 j_2} z^{2 j_3} ]
  =
  \int_{\theta = -\pi}^{\pi}
    \cos^{2 j_1} ( \theta )
    \sin^{2 j_2} ( \theta )
  \;
  d\theta
  \; \;
  \int_{\phi = -\pi/2}^{\pi/2}
    \cos^{2 j_1 + 2 j_2 + 1} ( \phi )
    \sin^{2 j_3} ( \phi )
  \;
  d\phi.
\end{eqnarray*}
These moments are calculated by \texttt{U3Moments.M}, a
\textsc{Mathematica} program (Appendix \ref{app:MathematicaCode}).
Table \ref{tab:U3Moments} lists them, for $m \leqslant 10 $, and this
allows us to establish $(*)$ for Gau{\ss} rules of degrees $ 3, 5,
\dots, 21 $.

\begin{table}[htbp]
  \centering
  \begin{tabular}{||c|c|c||c||c|c|c||c||}
    \hline\hline
    $ j_1 $  & $ j_2 $ & $ j_3 $ &
    $ \rule{0pt}{15pt} I [ x^{2 j_1} y^{2 j_2} z^{2 j_3} ] $ &
    $ j_1 $  & $ j_2 $ & $ j_3 $ &
    $ \rule{0pt}{15pt} I [ x^{2 j_1} y^{2 j_2} z^{2 j_3} ] $ \\
    \hline\hline
    $ 0 $&$ 0 $&$ ~1 $&$  4 \pi/3     $&$ 0 $&$ 4 $&$ ~6 $&$  4 \pi/12597 $ \\
    $ 0 $&$ 0 $&$ ~2 $&$  4 \pi/5     $&$ 0 $&$ 5 $&$ ~5 $&$ 12 \pi/46189 $ \\
    $ 0 $&$ 0 $&$ ~3 $&$  4 \pi/7     $&$ 1 $&$ 1 $&$ ~1 $&$  4 \pi/105 $ \\
    $ 0 $&$ 0 $&$ ~4 $&$  4 \pi/9     $&$ 1 $&$ 1 $&$ ~2 $&$  4 \pi/315 $ \\
    $ 0 $&$ 0 $&$ ~5 $&$  4 \pi/11    $&$ 1 $&$ 1 $&$ ~3 $&$  4 \pi/693 $ \\
    $ 0 $&$ 0 $&$ ~6 $&$  4 \pi/13    $&$ 1 $&$ 1 $&$ ~4 $&$  4 \pi/1287 $ \\
    $ 0 $&$ 0 $&$ ~7 $&$  4 \pi/15    $&$ 1 $&$ 1 $&$ ~5 $&$  4 \pi/2145 $ \\
    $ 0 $&$ 0 $&$ ~8 $&$  4 \pi/17    $&$ 1 $&$ 1 $&$ ~6 $&$  4 \pi/3315 $ \\
    $ 0 $&$ 0 $&$ ~9 $&$  4 \pi/19    $&$ 1 $&$ 1 $&$ ~7 $&$  4 \pi/4845 $ \\
    $ 0 $&$ 0 $&$ 10 $&$  4 \pi/21    $&$ 1 $&$ 1 $&$ ~8 $&$  4 \pi/6783 $ \\
    $ 0 $&$ 1 $&$ ~1 $&$  4 \pi/15    $&$ 1 $&$ 2 $&$ ~2 $&$  4 \pi/1155 $ \\
    $ 0 $&$ 1 $&$ ~2 $&$  4 \pi/35    $&$ 1 $&$ 2 $&$ ~3 $&$  4 \pi/3003 $ \\
    $ 0 $&$ 1 $&$ ~3 $&$  4 \pi/63    $&$ 1 $&$ 2 $&$ ~4 $&$  4 \pi/6435 $ \\
    $ 0 $&$ 1 $&$ ~4 $&$  4 \pi/99    $&$ 1 $&$ 2 $&$ ~5 $&$  4 \pi/12155 $ \\
    $ 0 $&$ 1 $&$ ~5 $&$  4 \pi/143   $&$ 1 $&$ 2 $&$ ~6 $&$  4 \pi/20995 $ \\
    $ 0 $&$ 1 $&$ ~6 $&$  4 \pi/195   $&$ 1 $&$ 2 $&$ ~7 $&$  4 \pi/33915 $ \\
    $ 0 $&$ 1 $&$ ~7 $&$  4 \pi/255   $&$ 1 $&$ 3 $&$ ~3 $&$  4 \pi/9009 $ \\
    $ 0 $&$ 1 $&$ ~8 $&$  4 \pi/323   $&$ 1 $&$ 3 $&$ ~4 $&$  4 \pi/21879 $ \\
    $ 0 $&$ 1 $&$ ~9 $&$  4 \pi/399   $&$ 1 $&$ 3 $&$ ~5 $&$  4 \pi/46189 $ \\
    $ 0 $&$ 2 $&$ ~2 $&$  4 \pi/105   $&$ 1 $&$ 3 $&$ ~6 $&$  4 \pi/88179 $ \\
    $ 0 $&$ 2 $&$ ~3 $&$  4 \pi/231   $&$ 1 $&$ 4 $&$ ~4 $&$ 28 \pi/415701 $ \\
    $ 0 $&$ 2 $&$ ~4 $&$  4 \pi/429   $&$ 1 $&$ 4 $&$ ~5 $&$  4 \pi/138567 $ \\
    $ 0 $&$ 2 $&$ ~5 $&$  4 \pi/715   $&$ 2 $&$ 2 $&$ ~2 $&$  4 \pi/5005 $ \\
    $ 0 $&$ 2 $&$ ~6 $&$  4 \pi/1105  $&$ 2 $&$ 2 $&$ ~3 $&$  4 \pi/15015 $ \\
    $ 0 $&$ 2 $&$ ~7 $&$  4 \pi/1615  $&$ 2 $&$ 2 $&$ ~4 $&$  4 \pi/36465 $ \\
    $ 0 $&$ 2 $&$ ~8 $&$  4 \pi/2261  $&$ 2 $&$ 2 $&$ ~5 $&$ 12 \pi/230945 $ \\
    $ 0 $&$ 3 $&$ ~3 $&$ 20 \pi/3003  $&$ 2 $&$ 2 $&$ ~6 $&$  4 \pi/146965 $ \\
    $ 0 $&$ 3 $&$ ~4 $&$  4 \pi/1287  $&$ 2 $&$ 3 $&$ ~3 $&$  4 \pi/51051 $ \\
    $ 0 $&$ 3 $&$ ~5 $&$  4 \pi/2431  $&$ 2 $&$ 3 $&$ ~4 $&$  4 \pi/138567 $ \\
    $ 0 $&$ 3 $&$ ~6 $&$  4 \pi/4199  $&$ 2 $&$ 3 $&$ ~5 $&$  4 \pi/323323 $ \\
    $ 0 $&$ 3 $&$ ~7 $&$  4 \pi/6783  $&$ 2 $&$ 4 $&$ ~4 $&$  4 \pi/415701 $ \\
    $ 0 $&$ 4 $&$ ~4 $&$ 28 \pi/21879 $&$ 3 $&$ 3 $&$ ~3 $&$ 20 \pi/969969 $ \\
    $ 0 $&$ 4 $&$ ~5 $&$ 28 \pi/46189 $&$ 3 $&$ 3 $&$ ~4 $&$ 20 \pi/2909907 $ \\
    \hline\hline
  \end{tabular}
  \caption{Moments of the first few monomials over $ U_3 $.}
  \label{tab:U3Moments}
\end{table}


\subsubsection{Solution of the System of Moment Equations $(*)$}

Having constructed tables of moments and possible structures for
various $ m $, we may consider solution of (the pared-down version of)
$( * ) $ for $ U_3 $.  Attempting to use \textsc{Mathematica} to do
this succeeded for $ m $ up to $ 5 $ (sometimes with a little human
assistance in making substitutions).  Beyond $ m = 5 $, the high degree
of the polynomials involved means that $(*)$ in general has only
transcendental solutions.  Approximate solution of $(*)$ is
attempted using numerical software in \textsc{matlab} (and \textsc{C}).

We use the \textsc{matlab} routine \texttt{fsolve}, which numerically
approximates the solution of a system of equations. This routine
requires a user-specified function that evaluates the system $( * ) $.
Initially, this was done as a \textsc{matlab} \texttt{m}-file.  As the
size of $(*)$ increases rapidly with $ m $, the large number of
expensive function evaluations involved made this a slow procedure.
Instead, \texttt{momenteq.c}, a \textsc{matlab} \texttt{mex}-file was
written (Appendix \ref{app:matlabCode}), which yielded a speedup of two
orders of magnitude. As numerous experiments are required to find a
satisfactory structure, this speedup is important.  To aid debugging, a
program \texttt{writestar.c} (Appendix \ref{app:CCode}), takes inputs
of $ m $ and $\left\lbrace K_i \right\rbrace $, and outputs a \LaTeX \,
file containing $(*)$ as a set of displayed equations.

A \textsc{matlab} driver program, \texttt{cubature.m} (Appendix
\ref{app:matlabCode}) is used to try various structures and choices of
$ m $.  For each $m $, \texttt{cubature.m} is run until either
\texttt{fsolve} solves $( * ) $, or the user gives up in disgust.
Successful results for $m = 1, \dots, 8 $ are presented in Appendix
\ref{app:GaussRulesU3}.

\pagebreak


\section{Product Rules for $ U_3 $}
\label{sec:ProductRules}

Product rules are an important approach to multidimensional cubature,
and they are the main direct competitor with Gau{\ss} rules. Here we
describe their construction for $ U_3 $, and compare their efficiency
with the Gau{\ss} rules discussed in \S\ref{sec:U3Cubature}. The material
is largely abstracted from chapter 2 of Stroud \cite{Stroud:71}.

\subsection{Introduction}
\label{sec:ProductRulesIntro}

To introduce the concept of product rules, consider the case where
$\mathcal{R}_n $ is $ C_3 $ (the unit cube), and we have a unit weighting
function. Let $ {\left\lbrace x_i, w_i \right\rbrace}_{i=1}^m $ be the
$m $-point 1D Gau{\ss}-Legendre rule of degree $ 2m-1 $, which exactly
integrates all of:
\begin{eqnarray*}
  \int_{-1}^1
    x^j
  dx
  \qquad
  j = 0, \dots, 2m - 1.
\end{eqnarray*}
For $ C_3 $, we wish to integrate exactly:
\begin{eqnarray*}
  \int_{-1}^1 \int_{-1}^1 \int_{-1}^1
    x_1^{j_1}
    x_2^{j_2}
    x_3^{j_3}
  \;
  dx_1 dx_2 dx_3.
\end{eqnarray*}
Writing this as an iterated integral allows us to construct a
\emph{product rule}, which exactly integrates all of:
\begin{eqnarray*}
  \int_{-1}^1
    x_1^{j_1} \;
  dx_1 \;
  \int_{-1}^1
    x_2^{j_2} \;
  dx_2 \;
  \int_{-1}^1
    x_3^{j_3} \;
  dx_3
  \qquad \qquad
  j_1, j_2, j_3 = 0, \dots, 2m - 1.
\end{eqnarray*}
Let $ \mathbf{i} = ( i_1, i_2, i_3 ) $ for all
$1 \leqslant i_1, i_2, i_3 \leqslant m $. We construct an
$m^3 $-point rule
$
  \left\lbrace {\mathbf{x}}_{\mathbf{i}}, A_{\mathbf{i}} \right\rbrace
$
for $ C_3 $ via:
\begin{eqnarray*}
  {\mathbf{x}}_{\mathbf{i}}
  =
  {(
    x_{i_1},
    x_{i_2},
    x_{i_3}
  )}^{\top}
  \qquad
  \mbox{and}
  \qquad
  A_{\mathbf{i}}
  =
  w_{i_1}
  w_{i_2}
  w_{i_3}.
\end{eqnarray*}
The set of points is a Cartesian product:
\begin{eqnarray*}
  \left\lbrace {\mathbf{x}}_{\mathbf{i}} \right\rbrace
  =
  \left\lbrace x_{i_1} \right\rbrace
  \times
  \left\lbrace x_{i_2} \right\rbrace
  \times
  \left\lbrace x_{i_3} \right\rbrace.
\end{eqnarray*}
This product rule will integrate exactly all monomials:
$
  x_1^{j_1} x_2^{j_2} x_3^{j_3}
$, for
$
  j_1, j_2, j_3 = 0, \dots, 2m - 1
$.
The highest degree monomial that it will integrate exactly is\\
\noindent$x_1^{2m-1} x_2^{2m-1} x_3^{2m-1} $, of degree $ 3 ( 2m - 1 )
$.  It is \emph{not} however, a rule of degree $ 3 ( 2m - 1 ) $ as it
does not exactly integrate \emph{all\/} polynomials of this degree, for
instance, it does not exactly integrate $ x_1^{3 ( 2m - 1 )} $. It
\emph{is} a Gau{\ss} rule in that it exactly integrates all polynomials
of up to a certain degree, but only in one variable.  However, as the
1D Gau{\ss} rules are good, this product rule is also good. This is a
general property of product rules.

The rule requires $ m^3 $ points. In general, a product rule on
$\mathcal{R}_n $, of degree $ 2m - 1 $ in each of the $ n $ variables,
will require $ N = m^n $ points.  This exponential growth in $ N $ was
called the `curse of dimensionality' by authors in the 1960s, when it
was believed that there was no escape from it.

\enlargethispage{\baselineskip}

\begin{quotation}
  For small $ n $ product formulas are very useful. For example
  if one wanted a subroutine that used a fixed $ 1000 $-point
  formula for a wide class of integrands for the $ 3 $-cube,
  $ w ( x, y, z ) = 1 $, we believe that one could do no better
  than the product of three copies of the $ 10 $-point
  Gau{\ss}-Legendre formula. This formula has degree $ 19 $ and
  there are, in fact, no nineteenth-degree formulas known for the
  $ 3 $-cube using fewer than $ 1000 $ points. (\dots a lower
  bound for the number of points in such a formula is $ 221 $.)
  \begin{flushright}
    Stroud (1971), \cite[p 25]{Stroud:71}.
  \end{flushright}
\end{quotation}

\pagebreak

In fact \cite{MantelRabinowitz:77} shows that in theory, an FSM rule
for $ C_3 $ can be constructed using $ 345 $ points, although there is
no explicit calculation of such. It seems likely that an FSGM rule on
less than $ 400 $ points exists.  This book \cite{Stroud:71}
illustrates the limitations on the usefulness of product rules in terms
of the relationship $ N = m^n $.  Table \ref{tab:ProductUsefulness}
illustrates some upper limits on $ m $ for various $ n $ if a ceiling
of $ N $ function evaluations is enforced.

\begin{table}
  \centering
  \begin{tabular}{||c||*{3}{c|c||}}
    \hline\hline
    $ n $ & $ m $ & $ N = m^n $ & $ m $ & $ N = m^n $ & $ m $ &
    $ N = m^n $ \\
    \hline\hline
    $ ~2 $ & $ 31 $ & $ 961 $ & $ 999 $ & $ 998001 $ &
      $ 31622 $ & $ 999950884 $ \\
    $ ~3 $ & $ ~9 $ & $ 729 $ & $ ~99 $ & $ 970299 $ &
      $ ~~999 $ & $ 997002999 $ \\
    $ ~4 $ & $ ~5 $ & $ 625 $ & $ ~31 $ & $ 923521 $ &
      $ ~~177 $ & $ 981506241 $ \\
    $ ~5 $ & $ ~3 $ & $ 243 $ & $ ~15 $ & $ 759375 $ &
      $ ~~~63 $ & $ 992436543 $ \\
    $ ~6 $ & $ ~3 $ & $ 729 $ & $ ~~9 $ & $ 531441 $ &
      $ ~~~31 $ & $ 887503681 $ \\
    $ ~7 $ & $ ~2 $ & $ 128 $ & $ ~~7 $ & $ 823543 $ &
      $ ~~~19 $ & $ 893871739 $ \\
    $ ~8 $ & $ ~2 $ & $ 256 $ & $ ~~5 $ & $ 390625 $ &
      $ ~~~13 $ & $ 815730721 $ \\
    $ ~9 $ & $ ~2 $ & $ 512 $ & $ ~~4 $ & $ 262144 $ &
      $ ~~~~9 $ & $ 387420489 $ \\
    $ 10 $ & $ ~1 $ & $ ~~1 $ & $ ~~3 $ & $ ~59049 $ &
      $ ~~~~7 $ & $ 282475249 $ \\
    $ 11 $ & $ ~1 $ & $ ~~1 $ & $ ~~3 $ & $ 177147 $ &
      $ ~~~~6 $ & $ 362797056 $ \\
    $ 12 $ & $ ~1 $ & $ ~~1 $ & $ ~~3 $ & $ 531441 $ &
      $ ~~~~5 $ & $ 244140625 $ \\
    $ 13 $ & $ ~1 $ & $ ~~1 $ & $ ~~2 $ & $ ~~8192 $ &
      $ ~~~~4 $ & $ ~67108864 $ \\
    $ 14 $ & $ ~1 $ & $ ~~1 $ & $ ~~2 $ & $ ~16384 $ &
      $ ~~~~4 $ & $ 268435456 $ \\
    $ 15 $ & $ ~1 $ & $ ~~1 $ & $ ~~2 $ & $ ~32768 $ &
      $ ~~~~3 $ & $ ~14348907 $ \\
    $ 16 $ & $ ~1 $ & $ ~~1 $ & $ ~~2 $ & $ ~65536 $ &
      $ ~~~~3 $ & $ ~43046721 $ \\
    $ 17 $ & $ ~1 $ & $ ~~1 $ & $ ~~2 $ & $ 131072 $ &
      $ ~~~~3 $ & $ 129140163 $ \\
    $ 18 $ & $ ~1 $ & $ ~~1 $ & $ ~~2 $ & $ 262144 $ &
      $ ~~~~3 $ & $ 387420489 $ \\
    $ 19 $ & $ ~1 $ & $ ~~1 $ & $ ~~2 $ & $ 524288 $ &
      $ ~~~~2 $ & $ ~~~524288 $ \\
    $ 20 $ & $ ~1 $ & $ ~~1 $ & $ ~~1 $ & $ ~~~~~1 $ &
      $ ~~~~2 $ & $ ~~1048576 $ \\
    \hline\hline
  \end{tabular}
  \caption
  {
    Limits on $ m $ for various $ n $ using product rules requiring
    $ N = m^n $ evaluation points, if ceilings of
    $ N < {10}^3, {10}^6 $ and $ {10}^9 $ are enforced.  (Modelled
    after table 2.1 in \protect\cite[p~24]{Stroud:71}.)
  }
  \label{tab:ProductUsefulness}
\end{table}

\pagebreak


\subsection{Product Rules for $ U_3 $}

Any point in $ U_3 $ can be uniquely characterised by the longitude
$\theta $, and the co-latitude $ \phi $. This simplifies the
construction of product rules, as we can express such a rule as the
product of rules that integrate over $ \theta $ and $ \phi $,
respectively.  The following construction is abstracted from
\cite[pp~34-35~and~40-41]{Stroud:71}.

For $ k = 1, 2 $, let $ \left\lbrace y_{k, i}, A_{k, i} \right\rbrace $ be the points and
weights in the $ m $-point 1D Gau{\ss}-Jacobi rules:
\begin{eqnarray*}
  \int_{-1}^1
    {(
      1 - y_k^2
    )}^{(k-2)/2}
    f ( y_k )
  d y_k
  \approx
  \sum_{i=1}^m
    A_{k, i}
    f ( y_{k, i} ).
\end{eqnarray*}
For $ k = 1 $, this is the Gau{\ss}-Chebyshev rule of the first
kind $ \left\lbrace y_{1,i}, A_{1,i} \right\rbrace $:
\begin{eqnarray*}
  \int_{-1}^1
    {(
      1 - y_1^2
    )}^{-1/2}
    f ( y_1 )
  d y_1
  \approx
  \sum_{i=1}^m
    A_{1,i}
    f ( y_{1,i} ).
\end{eqnarray*}
For $ i = 1, \dots, m $, these are \cite[p~114]{Krylov:62}, given by
the following formula (note that the weights $ A_{1,i} $ are constant):
\begin{eqnarray*}
  y_{1,i}
  =
  \cos
  (
    \frac{( 2i - 1 ) \pi}{2m}
  )
  \quad
  \mbox{and}
  \quad
  A_{1,i}
  =
  \pi / m.
\end{eqnarray*}

For $ k = 2 $, the formula is a Gau{\ss}-Legendre rule
$\left\lbrace y_{2,i}, A_{2,i} \right\rbrace $:
\begin{eqnarray*}
  \int_{-1}^1
    f ( y_2 )
  d y_2
  \approx
  \sum_{i=1}^m
    A_{2,i}
    f ( y_{2,i} ).
\end{eqnarray*}
This is not expressible in a simple closed form, but calculation is
routine and efficient (e.g. the implementation in \texttt{gauss.m} in
Appendix \ref{app:matlabCode}).

Now let $ \mathbf{i} \equiv ( i_1, i_2 ) $, for
$1 \leqslant i_1, i_2 \leqslant m $, and define
${\nu}_{\mathbf{i},1}, {\nu}_{\mathbf{i},2}, {\nu}_{\mathbf{i},3} $ by:
\begin{eqnarray*}
  {\nu}_{\mathbf{i},1}
  & = &
  \pm
  {( 1 - y_{2,i_2}^2 )}^{1/2}
  {( 1 - y_{1,i_1}^2 )}^{1/2}
  \\
  {\nu}_{\mathbf{i},2}
  & = &
  \pm
  {( 1 - y_{2,i_2}^2 )}^{1/2}
  y_{1,i_1}
  \\
  {\nu}_{\mathbf{i},3}
  & = &
  \pm y_{2,i_2}.
\end{eqnarray*}
A $ 2 m^2 $-point product rule
$\left\lbrace {\mathbf{x}}_{\mathbf{i}}, B_{\mathbf{i}} \right\rbrace $
of degree $ 2m - 1 $ for $ U_3 $ is then given by:
\begin{eqnarray*}
  {\mathbf{x}}_{\mathbf{i}}
  =
  {(
    {\nu}_{\mathbf{i},1}, {\nu}_{\mathbf{i},2}, {\nu}_{\mathbf{i},3}
  )}^{\top}
  \quad
  \mbox{and}
  \quad
  B_{\mathbf{i}}
  =
  A_{1,i_1} A_{2,i_2}.
\end{eqnarray*}
Some substitutions and relabelling in the construction shows that it is
actually simple, and we write the computations as an algorithm.

\pagebreak

\subsubsection*{Algorithm to Construct Product Rules for $ U_3 $}

\begin{enumerate}
\item
  Given $ m $, create a matrix $ \mathbf{i} \equiv ( i_1, i_2 ) $, where
  $ 1 \leqslant i_1, i_2 \leqslant m $. This matrix labels indices of
  the points in the rule.
\item
  Compute the $ m $-point Gau{\ss}-Legendre rule
  $ \left\lbrace y_{i_2}, A_{i_2} \right\rbrace $.
\item
  Set
  \begin{eqnarray*}
    {\mathbf{x}}_{\mathbf{i}}
    =
    \left[
    \begin{array}{c}
      \pm
      {( 1 - y_{i_2}^2 )}^{1/2}
      \sin ( \frac{( 2 i_1 - 1 ) \pi}{2m} )
      \\
      \pm
      {( 1 - y_{i_2}^2 )}^{1/2}
      \cos ( \frac{( 2 i_1 - 1 ) \pi}{2m} )
      \\
      \pm y_{i_2}
    \end{array}
    \right]
    \qquad
    \mbox{and}
    \qquad
    B_{\mathbf{i}}
    =
    \frac{\pi}{m}
    A_{i_2},
  \end{eqnarray*}
  as the $ 2 m^2 $ evaluation points, and the corresponding weights
  of a product rule of degree $ 2m - 1 $.
\end{enumerate}

Implementation as two \textsc{matlab} files is presented in
Appendix \ref{app:matlabCode}. Specific illustration of the points and
weights of the rules generated is rather pointless, but a table of
errors for integrating the appropriate polynomials generated by the
function \texttt{u3prod.m} shows it to work perfectly (the errors are
of order machine precision).

\pagebreak


\subsection{Comparison with Gau{\ss} Rules for $ U_3 $}

For the the case of $ U_3 $, the curse of dimensionality for product
rules is not the terrible scourge that it might have been.
Table \ref{tab:ComparisonU3} compares the number of points required with
varying degree for FSGM (or even FSMG) rules, as compared with the
product rules. In many applications, we may have no reason to want
rules of degree more than about $ 10 $.  Observe that there is only a
small factor of inefficiency in using the product rules, not the many
orders of magnitude that appear when the dimension is higher.

\begin{table}[htbp]
  \centering
  \begin{tabular}{||c||c|c|c||}
    \hline\hline
           & \multicolumn{2}{c|}{$ N $} & Ratio \\
    Degree & FSMG & Product             & \% \\
    \hline\hline
    $ ~3 $ & $ ~~6 $ & $ ~~8 $ & $ 75 $ \\
    $ ~5 $ & $ ~14 $ & $ ~18 $ & $ 78 $ \\
    $ ~7 $ & $ ~26 $ & $ ~32 $ & $ 81 $ \\
    $ ~9 $ & $ ~38 $ & $ ~50 $ & $ 76 $ \\
    $ 11 $ & $ ~50 $ & $ ~72 $ & $ 69 $ \\
    $ 13 $ & $ ~78 $ & $ ~98 $ & $ 80 $ \\
    $ 15 $ & $ ~86 $ & $ 128 $ & $ 67 $ \\
    $ 17 $ & $ 110 $ & $ 162 $ & $ 68 $ \\
    \hline\hline
  \end{tabular}
  \caption
  {
    Comparison of the number of points required for the FSGM (or FSMG)
    Gau{\ss} rules and the product rules for $ U_3 $. The last column
    is the ratio of the number of points (function evaluations)
    required by the Gau{\ss} rules relative to the product rules.
  }
  \label{tab:ComparisonU3}
\end{table}


\subsection{Alternative Philosophies}
\label{sec:Alternatives}

Gau{\ss} rules (and their ancestors the equally-spaced formulae) are
based on exploiting the analytical properties of smooth integrands.
They are optimal in the sense that in general they will be the best
choice for approximating integrals involving smooth integrands.  As
mentioned in \S\ref{sec:FullSymmetry}, the optimal multidimensional
cubature can only be sensibly considered for regions with some
symmetry, and we dealt with the case of full symmetry.

To deal with non-smooth integrands, possibly even random distributions,
and with non-symmetric, possibly even disconnected domains of
integration, alternative philosophies are usually more relevant.
Textbooks on numerical integration commonly contain many pages
describing minute implementational details to further refine the theory
for optimal methods; and then devote a similar amount of space to
describing real alternatives (e.g.
\cite{DavisRabinowitz:84,Stroud:71}).  The main alternative technique
is the Monte Carlo method, but a more recent idea is the `lattice'
method.


\subsubsection{Monte Carlo `Simulation' Methods}

These techniques are based on averaging function values at a random
selection of points within the region, and they are particularly
appropriate for oddly-shaped regions and non-smooth integrands. They
are widely used in statistical applications, where integrals of high
dimension must be approximated. Performance is often about $ N^{-1} $
(the error incurred using $ N $ points should be proportional to
$N^{-1} $). Volumes have been written about them (e.g. see
\cite[chapter~6]{Stroud:71}).


\subsubsection{Lattice Methods}

These methods generalise the idea of placement of equally-spaced points
in 1D (with weights selected according to some generally simple
formula, expressible in closed form). The idea is to catch as
representative a sample of function values as possible. In 1D, this leads
to the Newton-Cotes family of rules.  For smooth integrands, these
rules increase in accuracy algebraically with their number of
evaluation points, although they are not always good.  (Whilst this
may be adequate, it is still inferior to the exponential accuracy of
the Gau{\ss} rules.)

The problem with attempting to generalise this to the case of several
dimensions is that the notion of `equal spacing' of evaluation points
becomes less well-defined.  Placement of equally-spaced points on
$U_3 $ is equivalent to maximally covering it with nonintersecting
equal circles, a problem thought to be intractable.  The `best' that
can be done involves heuristic algorithms, and lots of computer time
\cite
{%
  ClareKeppert:86,%
  ClareKeppert:91%
}.

Nevertheless, it may be reasonable to try to approximate the
equally-spaced placement of points within our region.  Recent research
involving Sloan and Lyness
\cite
{%
  LynessSloan:89,%
  Sloan:91,%
  SloanLyness:89%
}, has achieved this using geometric construction techniques, called
`lattice methods'.  For periodic functions on $ {[ 0, 1 ]}^N $,
lattice methods generalise the trapezoidal rule, preserving the order
of the error as $ N^{-2} $.  For $ U_3 $, a placement called `spherical
$t $-designs' is used \cite{Sloan:91}.


\subsubsection{Comparison with Gau{\ss} Cubature}

For FS regions and weight functions, where the expected integrands are
smooth, these alternatives are a poor second choice in comparison with
Gau{\ss} rules. The Monte Carlo methods will only converge as $ N^{-1}
$, and the theory for the lattice methods is not very general, results
only being available for one type of region at a time.

It cannot be overemphasised that where Gau{\ss} cubature
is available, it should be used, particularly as the dimensionality
increases.  For evaluation of integrals over 2D manifolds, Gau{\ss}
cubature is applicable.

\pagebreak


\appendix

\section{Number of Equivalence Classes with Dimension}
\label{app:ECData}

\begin{figure}[htbp]
  \begin{center}
    \epsfig{file=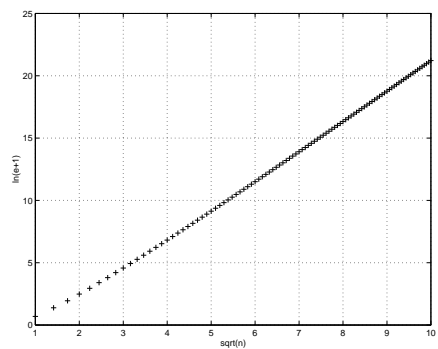,height=130pt}
    \caption{Plot of the data in Table \protect\ref{tab:NumECTypes}.}
    \label{fig:ECData}
  \end{center}
\end{figure}

Results from running \texttt{findec.c} (Appendix \ref{app:CCode}) for
$n = 1, \dots, 100 $ are listed in Table \ref{tab:NumECTypes}, and
graphically presented in Figure \ref{fig:ECData}.  For $ n = 100 $, the
(optimised) program requires about $ 1000 $ CPU minutes on a
\textsc{SPARC-10} workstation.  Data has been manually checked for
$n = 1, \dots, 10 $.

\begin{table}[htbp]
  \centering
  \begin{tabular}{||*{5}{*{2}{c|}|}}
    \hline\hline
    $ n $ & $ e+1 $ & $ n $ & $ e+1 $ & $ n $ & $ e+1 $ &
    $ n $ & $ e+1 $ & $ n $ & $ e+1 $ \\
    \hline\hline
$~1 $&$ ~~~2 $&$ 21 $&$ ~~3506 $&$ 41 $&$ ~259891 $&$ 61 $&$ ~~7760854 $&$ ~81 $&$ ~141227966 $ \\
$~2 $&$ ~~~4 $&$ 22 $&$ ~~4508 $&$ 42 $&$ ~313065 $&$ 62 $&$ ~~9061010 $&$ ~82 $&$ ~161734221 $ \\
$~3 $&$ ~~~7 $&$ 23 $&$ ~~5763 $&$ 43 $&$ ~376326 $&$ 63 $&$ ~10566509 $&$ ~83 $&$ ~185072690 $ \\
$~4 $&$ ~~12 $&$ 24 $&$ ~~7338 $&$ 44 $&$ ~451501 $&$ 64 $&$ ~12308139 $&$ ~84 $&$ ~211616350 $ \\
$~5 $&$ ~~19 $&$ 25 $&$ ~~9296 $&$ 45 $&$ ~540635 $&$ 65 $&$ ~14320697 $&$ ~85 $&$ ~241783707 $ \\
$~6 $&$ ~~30 $&$ 26 $&$ ~11732 $&$ 46 $&$ ~646193 $&$ 66 $&$ ~16644217 $&$ ~86 $&$ ~276046669 $ \\
$~7 $&$ ~~45 $&$ 27 $&$ ~14742 $&$ 47 $&$ ~770947 $&$ 67 $&$ ~19323906 $&$ ~87 $&$ ~314934342 $ \\
$~8 $&$ ~~67 $&$ 28 $&$ ~18460 $&$ 48 $&$ ~918220 $&$ 68 $&$ ~22411641 $&$ ~88 $&$ ~359042451 $ \\
$~9 $&$ ~~97 $&$ 29 $&$ ~23025 $&$ 49 $&$ 1091745 $&$ 69 $&$ ~25965986 $&$ ~89 $&$ ~409038376 $ \\
$10 $&$ ~139 $&$ 30 $&$ ~28629 $&$ 50 $&$ 1295971 $&$ 70 $&$ ~30053954 $&$ ~90 $&$ ~465672549 $ \\
$11 $&$ ~195 $&$ 31 $&$ ~35471 $&$ 51 $&$ 1535914 $&$ 71 $&$ ~34751159 $&$ ~91 $&$ ~529784908 $ \\
$12 $&$ ~272 $&$ 32 $&$ ~43820 $&$ 52 $&$ 1817503 $&$ 72 $&$ ~40143942 $&$ ~92 $&$ ~602318715 $ \\
$13 $&$ ~373 $&$ 33 $&$ ~53963 $&$ 53 $&$ 2147434 $&$ 73 $&$ ~46329631 $&$ ~93 $&$ ~684328892 $ \\
$14 $&$ ~508 $&$ 34 $&$ ~66273 $&$ 54 $&$ 2533589 $&$ 74 $&$ ~53419131 $&$ ~94 $&$ ~776998612 $ \\
$15 $&$ ~684 $&$ 35 $&$ ~81156 $&$ 55 $&$ 2984865 $&$ 75 $&$ ~61537395 $&$ ~95 $&$ ~881650031 $ \\
$16 $&$ ~915 $&$ 36 $&$ ~99133 $&$ 56 $&$ 3511688 $&$ 76 $&$ ~70826486 $&$ ~96 $&$ ~999764335 $ \\
$17 $&$ 1212 $&$ 37 $&$ 120770 $&$ 57 $&$ 4125842 $&$ 77 $&$ ~81446349 $&$ ~97 $&$ 1132995265 $ \\
$18 $&$ 1597 $&$ 38 $&$ 146785 $&$ 58 $&$ 4841062 $&$ 78 $&$ ~93578513 $&$ ~98 $&$ 1283193401 $ \\
$19 $&$ 2087 $&$ 39 $&$ 177970 $&$ 59 $&$ 5672882 $&$ 79 $&$ 107427163 $&$ ~99 $&$ 1452423276 $ \\
$20 $&$ 2714 $&$ 40 $&$ 215308 $&$ 60 $&$ 6639349 $&$ 80 $&$ 123223639 $&$ 100 $&$ 1642992568 $ \\
    \hline\hline
  \end{tabular}
  \caption
  {
    The number $ e + 1 $ of types of equivalence classes of FS sets of
    points for $ n = 1, \dots, 100 $.
  }
  \label{tab:NumECTypes}
\end{table}

\vspace{5mm}

\clearpage


\section{Gau{\ss} Cubature Rules for $ U_3 $}
\label{app:GaussRulesU3}

We present computed FS (usually FSMG) Gau{\ss} rules for $ U_3 $, for
$m = 1, \dots, 8 $ (degrees $ 3, 5, \dots, 17 $). For $ m \leqslant 5 $,
\textsc{Mathematica} (or a combination of it and some manual
substitutions) provides analytic solutions.  For higher degrees, only
transcendental solutions exist. Instead, we use the \textsc{matlab}
program \texttt{cubature.m} to approximate a solution, and we have to
trust that an approximation with a small residual corresponds to a
transcendental solution.  The analytic solutions for low degrees
provide good test data for \texttt{cubature.m}.  All of the rules
presented have been discovered by \texttt{cubature.m}; where possible,
analytic solutions have been substituted. Beyond $ m = 8 $, no
solutions at all have been found, but they should be available with
sufficient computational effort.  For most cases the rules are FSMG
from the first structure corresponding to the first minima of the IPP.
For $ m = 4 $, the structure is from the second minima of the first
structure. For $ m = 6 $ there is no FSMG rule, but an FSGM rule from
the first minima of the second structure is obtained.


\begin{itemize}
\item
  $ m = 1 $, degree $ 3 $. An FSMG rule is
  $
    U_3 : 3 \mbox{--} 1.1 ( 1, 0, 0, 0, 0, 0 ) \mbox{--} 6
  $:
  \begin{eqnarray*}
    a_1 & = & \frac{2 \pi}{3}
    \qquad \qquad
    \alpha_1 = 1.
  \end{eqnarray*}
\item
  $ m = 2 $, degree $ 5 $. An FSMG rule is
  $
    U_3 : 5 \mbox{--} 1.1 ( 1, 0, 0, 1, 0, 0 ) \mbox{--} 14
  $:
  \begin{eqnarray*}
    a_1 & = & \frac{4 \pi}{15}
    \qquad \qquad
    \alpha_1 = 1
    \\
    d_1 & = & \frac{3 \pi}{10}
    \qquad \qquad
    \epsilon_1 = \frac{1}{\sqrt{3}}.
  \end{eqnarray*}
\item
  $ m = 3 $, degree $ 7 $. An FSMG rule is
  $
    U_3 : 7 \mbox{--} 1.1 ( 1, 1, 0, 1, 0, 0 ) \mbox{--} 26
  $:
  \begin{eqnarray*}
    a_1 & = & \frac{4 \pi}{21}
    \qquad \qquad
    \alpha_1 = 1
    \\
    b_1 & = & \frac{16 \pi}{105}
    \qquad \quad \;
    \beta_1 = \frac{1}{\sqrt{2}}
    \\
    d_1 & = & \frac{9 \pi}{70}
    \qquad \qquad
    \epsilon_1 = \frac{1}{\sqrt{3}}.
  \end{eqnarray*}
\item
  $ m = 4 $, degree $ 9 $. An FSMG rule is
  $
    U_3 : 9 \mbox{--} 1.2 ( 1, 0, 1, 1, 0, 0 ) \mbox{--} 38
  $:
  \begin{eqnarray*}
    a_1 & = & \frac{4 \pi}{105}
    \qquad \qquad
    \alpha_1 = 1
    \\
    c_1 & = & \frac{4 \pi}{35}
    \qquad \qquad \;
    \gamma_1 = \sqrt{ \frac{1}{2} ( 1 - \frac{1}{\sqrt{3}} )}
    \qquad \qquad
    \delta_1 = \sqrt{ \frac{1}{2} ( 1 + \frac{1}{\sqrt{3}} )}
    \\
    d_1 & = & \frac{9 \pi}{70}
    \qquad \qquad \;
    \epsilon_1 = \frac{1}{\sqrt{3}}.
  \end{eqnarray*}

\pagebreak

\item
  $ m = 5 $, degree $ 11 $. An FSMG rule is
  $
    U_3 : 11 \mbox{--} 1.1 ( 1, 1, 0, 1, 1, 0 ) \mbox{--} 50
  $:
  \begin{eqnarray*}
    a_1 & = & \frac{16 \pi}{315}
    \qquad \qquad \; \;
    \alpha_1 = 1
    \\
    b_1 & = & \frac{256 \pi}{2835}
    \qquad \qquad
    \beta_1 = \frac{1}{\sqrt{2}}
    \\
    d_1 & = & \frac{27 \pi}{320}
    \qquad \qquad \; \;
    \epsilon_1 = \frac{1}{\sqrt{3}}
    \\
    e_1 & = & \frac{14641 \pi}{181440}
    \qquad \quad
    \zeta_1 = \frac{1}{\sqrt{11}}
    \qquad \qquad
    \eta_1 = \frac{3}{\sqrt{11}}.
  \end{eqnarray*}
\item
  $ m = 6 $, degree $ 13 $. An FSM rule is
  $
    U_3 : 13 \mbox{--} 1.1 ( 1, 1, 1, 1, 1, 0 ) \mbox{--} 74
  $:
  \begin{eqnarray*}
    a_1 & \approx & 0.00644739233053
    \qquad \qquad
    \alpha_1 = 1
    \\
    b_1 & \approx & 0.20865289186971
    \qquad \qquad
    \beta_1 = 1 / \sqrt{2}
    \\
    c_1 & \approx & 0.20762372406088
    \qquad \qquad
    \gamma_1 \approx 0.32077264898077
    \qquad \qquad
    \delta_1 \approx 0.94715622136259
    \\
    d_1 & \approx & -0.37178913059595
    \qquad \quad
    \epsilon_1 = 1 / \sqrt{3}
    \\
    e_1 & \approx & 0.33396646771858
    \qquad \qquad
    \zeta_1 \approx 0.48038446141531
    \qquad \qquad
    \eta_1 \approx 0.73379938570528.
  \end{eqnarray*}
  This rule is not good, as $ d_1 < 0 $.  A (non-unique) FSGM rule is
  $
    U_3 : 13 \mbox{--} 2.1 ( 1, 0, 1, 0, 2, 0 ) \mbox{--} 78
  $:
  \begin{eqnarray*}
    a_1 & \approx & 0.05571838151106
    \qquad \qquad
    \alpha_1 = 1
    \\
    c_1 & \approx & 0.18861500631211
    \qquad \qquad
    \gamma_1 \approx 0.33370053800545
    \qquad \qquad
    \delta_1 \approx 0.94267913466612
    \\
    e_1 & \approx & 0.12537551702973
    \qquad \qquad
    \zeta_1 \approx 0.70117074174860
    \qquad \qquad
    \eta_1 \approx 0.12930267526790
    \\
    e_2 & \approx & 0.19567865687870
    \qquad \qquad
    \zeta_2 \approx 0.43948383947130
    \qquad \qquad
    \eta_2 \approx 0.78339511722191.
  \end{eqnarray*}
\item
  $ m = 7 $, degree $ 15 $. An FSMG rule is
  $
    U_3 : 15 \mbox{--} 1.1 ( 1, 0, 1, 1, 2, 0 ) \mbox{--} 86
  $:
  \begin{eqnarray*}
    a_1 & \approx & 0.14506632743849
    \qquad \qquad
    \alpha_1 = 1
    \\
    c_1 & \approx & 0.14843778669299
    \qquad \qquad
    \gamma_1 \approx 0.92733065715117
    \qquad \qquad
    \delta_1 \approx 0.37424303909034
    \\
    d_1 & \approx & 0.15009158815708
    \qquad \qquad
    \epsilon_1 = 1 / \sqrt{3}
    \\
    e_1 & \approx & 0.13961936079093
    \qquad \qquad
    \zeta_1 \approx 0.36960284645415
    \qquad \qquad
    \eta_1 \approx 0.85251831170127
    \\
    e_2 & \approx & 0.14924451686907
    \qquad \qquad
    \zeta_2 \approx 0.69435400660267
    \qquad \qquad
    \eta_2 \approx 0.18906355288540.
  \end{eqnarray*}
\item
  $ m = 8 $, degree $ 17 $. An FSMG rule is
  $
    U_3 : 17 \mbox{--} 1.1 ( 1, 0, 1, 1, 3, 0 ) \mbox{--} 110
  $:
  \begin{eqnarray*}
    a & \approx & 0.04810746585109
    \qquad \qquad
    \alpha = 1
    \\
    c_1 & \approx & 0.12183091738552
    \qquad \qquad
    \gamma_1 \approx 0.87815891060407
    \qquad \qquad
    \delta_1 \approx 0.47836902881214
    \\
    d_1 & \approx & 0.12307173528176
    \qquad \qquad
    \epsilon_1 = 1 / \sqrt{3}
    \\
    e_1 & \approx & 0.10319173408833
    \qquad \qquad
    \zeta_1 \approx 0.18511563534456
    \qquad \qquad
    \eta_1 \approx 0.96512403508666
    \\
    e_2 & \approx & 0.12058024902856
    \qquad \qquad
    \eta_2 \approx 0.82876998125269
    \qquad \qquad
    \zeta_2 \approx 0.39568947305584
    \\
    e_3 & \approx & 0.12494509687253
    \qquad \qquad
    \zeta_3 \approx 0.69042104838229
    \qquad \qquad
    \eta_3 \approx 0.21595729184587.
  \end{eqnarray*}
\end{itemize}

\pagebreak


\section{\textsc{C} Code}
\label{app:CCode}

\subsection{findec.c}
\input{C/findec}

\subsection{ipp.c}
\input{C/ipp}

\pagebreak

\subsection{writestar.c}
\input{C/writestar}

\vspace{10mm}

\input{C/eqnout}

\pagebreak


\section{\textsc{MATLAB Code}}
\label{app:matlabCode}

\subsection{$ U_3 $ Gau{\ss} Rules}

\subsubsection{cubature.m}
\input{matlab/cubature}

\subsubsection{momenteq.c}
\input{matlab/momenteq}


\subsection{$ U_3 $ Product Rules}

\subsubsection{u3prod.m}
\input{matlab/u3prod}

\pagebreak

\subsubsection{gauss.m}
\input{matlab/gauss}

\pagebreak


\section{{\textsc{Mathematica}} Code}
\label{app:MathematicaCode}

\subsection{IPPBasicSolution.M}
\input{Mathematica/IPPBasicSoln}

\pagebreak

\subsection{U3Moments.M}
\input{Mathematica/U3Moments}

\pagebreak


\bibliography{DeWit93}
\bibliographystyle{plain}

\end{document}

%% file: C/findec.tex
\small\tt 
\begin{verbatim}
/*
  Find the number of types of equivalence classes e+1, of fully
  symmetric points in n dimensions, by exhaustively enumerating them.
  A UNIX input line "findEC a b" will find e+1 for n = a, ..., b.

  Peter Adams and David  De Wit
  August 1  1993
*/

int string[100], ctr, n;

int checkrep(ind)
int ind;
{
  int repnums[100], rn, i, lv;

  for (i = 0; i < 100; i++)
    repnums[i] = 0;
  repnums[rn = 0] = 1;
  lv = 1;
  for (i = 1; i <= ind; i++) {
    if (string[i] == lv)
      repnums[rn]++;
    else {
      lv = string[i];
      repnums[++rn] = 1;
    }
  }
  for (i = 1; i < rn; i++)
    if (repnums[i] > repnums[i-1])
      return(0);
  return(1);
}

void build(lastv, ind)
int lastv, ind;
{
  int i;

  if (!checkrep(ind)) return;
  ctr++;
  if (ind == n) return;
  string[ind] = lastv;
  build(lastv, ind+1);
  if (lastv < n) {
    string[ind] = lastv + 1;
    build(lastv+1, ind+1);
  }
  string[ind] = 0;
}

main(argc, argv)
char *argv[];
{
  int i;
  
  printf(" n\t   e+1\n--------------\n", n, ctr);
  for (n = atoi(argv[1]); n <= atoi(argv[2]); n++)
  {
    ctr = 0;
    for (i = 0; i < n; i++)
      string[i] = 0;
    ctr++;
    string[0] = 1;
    build(1, 1);
    printf("%2d\t%6d\n", n, ctr);
  }
}
\end{verbatim} 
\normalsize\rm

%% file: C/ipp.tex
\small\tt 
\begin{verbatim}
/*
    Exhaustively solve the IPP, using the precalculated lower bound for
    $ N $. This is a crude but fast and successful method. The data is
    output, for each m, in terms of increasing structure. Running the
    UNIX "sort" on the output orders it into increasing N, then structure.

    David  De Wit
    August 1 -- August 9  1993
*/

#include <stdio.h>

#define LIMIT   20
#define NUMofM  21

static int NLB[NUMofM] =
{
  0,   6,  14,  26,  38,  50,  74,  86, 110, 138, 162,
     190, 230, 258, 298, 342, 382, 426, 482, 526, 582
};
static int c[4][NUMofM] =
{
  {0, 1, 2, 3, 4, 5, 7, 8, 10, 12, 14, 16, 19, 21, 24, 27, 30, 33, 37, 40, 44},
  {0, 0, 0, 1, 1, 2, 3, 4,  5,  7,  8, 10, 12, 14, 16, 19, 21, 24, 27, 30, 33},
  {0, 0, 0, 0, 0, 0, 1, 1,  2,  3,  4,  5,  7,  8, 10, 12, 14, 16, 19, 21, 24},
  {0, 0, 0, 0, 0, 0, 0, 0,  0,  1,  1,  2,  3,  4,  5,  7,  8, 10, 12, 14, 16}
};
static double alpha[NUMofM] =
{
  0.00, 17.00, 7.00, 4.00, 3.00, 2.30, 1.80, 1.70, 1.45, 1.40, 1.35,
         1.25, 1.20, 1.20, 1.15, 1.13, 1.12, 1.11, 1.09, 1.08, 1.08
};

main()
{
  int     m, K1, K2, K3, K4, K5, K6, N, Nvars, nstruct[NUMofM], NUB[NUMofM];

  for (m = 1; m < NUMofM; m++) {
    nstruct[m] = 0;
    NUB[m] = (int) (alpha[m]*NLB[m]);
    for (K1 = 0; K1 <= 1; K1++)
    for (K2 = 0; K2 <= 1; K2++)
    for (K3 = 0; K3 <= LIMIT; K3++)
    for (K4 = 0; K4 <= 1; K4++)
    for (K5 = 0; K5 <= LIMIT; K5++)
    for (K6 = 0; K6 <= LIMIT; K6++) {
      N = 6*K1 + 12*K2 + 24*K3 + 8*K4 + 24*K5 + 48*K6;
      if (N > NUB[m])
         break;
      if (NLB[m] <= N &
          K1 + K2 + 2*K3 + K4 + 2*K5 + 3*K6 >= c[0][m] &
                           K4 + 2*K5 + 3*K6 >= c[1][m] &
                    2*K3             + 3*K6 >= c[2][m] &
                                       3*K6 >= c[3][m]) {
        Nvars = 2*K1 + 2*K2 + 3*K3 + 2*K4 + 3*K5 + 4*K6;
        printf("$ %2d $ & $ %3d $ & $ %1d $ & $ %1d $ & $ ", m, N, K1, K2);
        printf("%2d $ & $ %1d $ & $ %2d $ & $ ", K3, K4, K5);
        printf("%2d $ & $ %2d $ \\\\\n", K6, Nvars);
        nstruct[m]++;
      }
    }
    printf("\\hline\n %% For m = %d, there are %d structures in the range %d-%d\n\n",
           m, nstruct[m], NLB[m], NUB[m]);
  }
}
\end{verbatim} 
\normalsize\rm

%% file: C/writestar.tex
\small\tt 
\begin{verbatim}
/*
  Write down the system of moment equations as a LaTeX file.  Input is
  m, and a structure K.  Includes the moment data created by the
  Mathematica function U3Moments.M. This code is not watertight!

  David  De Wit
  August 4  --  September 29  1993
*/

#include <stdio.h>
#define    pi    3.141592653589793238462643383280

main(argc, argv)
int  argc;
char *argv[];
{
  int    K[7], m, i, j1, j2, j3, J1, J2, J3, eqno, pflag, nflag;
  double mom[4][6][11];

/* Initialise some variables */
  if (argc == 8) {
    m = atoi(argv[1]);
    for (i = 1; i <= 6; i++)
      K[i] = atoi(argv[i+1]);
  }
  else {
    m = 3;  K[1] = 1;  K[2] = 1;  K[3] = 0;  K[4] = 1;  K[5] = 0;  K[6] = 0;
  }

  mom[0][0][0]  =  4*pi;            mom[0][0][1]  =  4*pi/3;
  mom[0][0][2]  =  4*pi/5;          mom[0][0][3]  =  4*pi/7;
  mom[0][0][4]  =  4*pi/9;          mom[0][0][5]  =  4*pi/11;
  mom[0][0][6]  =  4*pi/13;         mom[0][0][7]  =  4*pi/15;
  mom[0][0][8]  =  4*pi/17;         mom[0][0][9]  =  4*pi/19;
  mom[0][0][10] =  4*pi/21;         mom[0][1][1]  =  4*pi/15;
  mom[0][1][2]  =  4*pi/35;         mom[0][1][3]  =  4*pi/63;
  mom[0][1][4]  =  4*pi/99;         mom[0][1][5]  =  4*pi/143;
  mom[0][1][6]  =  4*pi/195;        mom[0][1][7]  =  4*pi/255;
  mom[0][1][8]  =  4*pi/323;        mom[0][1][9]  =  4*pi/399;
  mom[0][2][2]  =  4*pi/105;        mom[0][2][3]  =  4*pi/231;
  mom[0][2][4]  =  4*pi/429;        mom[0][2][5]  =  4*pi/715;
  mom[0][2][6]  =  4*pi/1105;       mom[0][2][7]  =  4*pi/1615;
  mom[0][2][8]  =  4*pi/2261;       mom[0][3][3]  =  20*pi/3003;
  mom[0][3][4]  =  4*pi/1287;       mom[0][3][5]  =  4*pi/2431;
  mom[0][3][6]  =  4*pi/4199;       mom[0][3][7]  =  4*pi/6783;
  mom[0][4][4]  =  28*pi/21879;     mom[0][4][5]  =  28*pi/46189;
  mom[0][4][6]  =  4*pi/12597;      mom[0][5][5]  =  12*pi/46189;
  mom[1][1][1]  =  4*pi/105;        mom[1][1][2]  =  4*pi/315;
  mom[1][1][3]  =  4*pi/693;        mom[1][1][4]  =  4*pi/1287;
  mom[1][1][5]  =  4*pi/2145;       mom[1][1][6]  =  4*pi/3315;
  mom[1][1][7]  =  4*pi/4845;       mom[1][1][8]  =  4*pi/6783;
  mom[1][2][2]  =  4*pi/1155;       mom[1][2][3]  =  4*pi/3003;
  mom[1][2][4]  =  4*pi/6435;       mom[1][2][5]  =  4*pi/12155;
  mom[1][2][6]  =  4*pi/20995;      mom[1][2][7]  =  4*pi/33915;
  mom[1][3][3]  =  4*pi/9009;       mom[1][3][4]  =  4*pi/21879;
  mom[1][3][5]  =  4*pi/46189;      mom[1][3][6]  =  4*pi/88179;
  mom[1][4][4]  =  28*pi/415701;    mom[1][4][5]  =  4*pi/138567;
  mom[2][2][2]  =  4*pi/5005;       mom[2][2][3]  =  4*pi/15015;
  mom[2][2][4]  =  4*pi/36465;      mom[2][2][5]  =  12*pi/230945;
  mom[2][2][6]  =  4*pi/146965;     mom[2][3][3]  =  4*pi/51051;
  mom[2][3][4]  =  4*pi/138567;     mom[2][3][5]  =  4*pi/323323;
  mom[2][4][4]  =  4*pi/415701;     mom[3][3][3]  =  20*pi/969969;
  mom[3][3][4]  =  20*pi/2909907;

/* Introduce the output */
  printf("Output from running `writestar': m = %d, ", m);
  printf("K = %d %d %d %d %d %d:\n\n", K[1], K[2], K[3], K[4], K[5], K[6]);

/* Deal with Subsystem I */
  printf("\\be\n\%% Subsystem I/1:\n");
  pflag = 0;  eqno++;
  printf("  I \\[ 1 \\]\n  \\eq\n  ");
  if (K[1]) {
    pflag = 1;  printf("6 a_1");
  }
  if (K[2]) {
    if (pflag) printf(" +\n  ");  pflag = 1;
    printf("12 b_1");
  }
  if (K[3]) {
    if (pflag) printf(" +\n  ");  pflag = 1;
    if (K[3] > 1)
      printf("24 \\sum_{i=1}^{%d} c_i", K[3]);
    else
      printf("24 c_1");
  }
  if (K[4]) {
    if (pflag) printf(" +\n  ");  pflag = 1;
    printf("8 d_1");
  }
  if (K[5]) {
    if (pflag) printf(" +\n  ");  pflag = 1;
    if (K[5] > 1)
      printf("24 \\sum_{i=1}^{%d} e_i", K[5]);
    else
      printf("24 e_1");
  }
  if (K[6]) {
    if (pflag) printf(" +\n  ");  pflag = 1;
    if (K[6] > 1)
      printf("48 \\sum_{i=1}^{%d} f_i", K[6]);
    else
      printf("48 f_1");
  }
  printf("\n  \\\\\n");

  printf("\%% Subsystem I/2:\n");
  for (j1 = 1; j1 <= m; j1++)
  {
    J1 = 2*j1;  eqno++;
    printf("  I \\[ x^{%d} \\]\n  \\eq\n  ", J1);
    if (K[1]) {
      pflag = 1;  printf("2 a_1 \\alpha_1^{%d}", J1);
    }
    if (K[2]) {
      if (pflag) printf(" +\n  ");  pflag = 1;
      printf("8 b_1 \\beta_1^{%d}", J1);
    }
    if (K[3]) {
      if (pflag) printf(" +\n  ");  pflag = 1;
      if (K[3] > 1)
        printf("8 \\sum_{i=1}^{%d} c_i \\( \\gamma_i^{%d} + \\delta_i^{%d} \\)", K[3], J1, J1);
      else
        printf("8 c_1 \\( \\gamma_1^{%d} + \\delta_1^{%d} \\)", J1, J1);
    }
    if (K[4]) {
      if (pflag) printf(" +\n  ");  pflag = 1;
      printf("8 d_1 \\epsilon_1^{%d}", J1);
    }
    if (K[5]) {
      if (pflag) printf(" +\n  ");  pflag = 1;
      if (K[5] > 1)
        printf("8 \\sum_{i=1}^{%d} e_i \\( 2 \\zeta_i^{%d} + \\eta_i^{%d} \\)", K[5], J1, J1);
      else
        printf("8 e_1 \\( 2 \\zeta_1^{%d} + \\eta_1^{%d} \\)", J1, J1);
    }
    if (K[6]) {
      if (pflag) printf(" +\n  ");  pflag = 1;
      if (K[6] > 1) {
        printf("16 \\sum_{i=1}^{%d} f_i ", K[6]);
        printf("\\( \\theta_i^{%d} + \\mu_i^{%d} + \\lambda_i^{%d} \\)", J1, J1, J1);
      }
      else
        printf("16 f_1 \\( \\theta_1^{%d} + \\mu_1^{%d} + \\lambda_1^{%d} \\)", J1, J1, J1);
    }
    printf("\n  \\\\\n");
  }

/* Deal with Subsystem II */
  if (m >= 2) {
  printf("\%% Subsystem II:\n");
  for (j1 = 1; j1 <= m; j1++)
  for (j2 = j1; j2 <= m-j1; j2++) {
    pflag = 0;  J1 = 2*j1;  J2 = 2*j2;
    eqno++;
    printf("  I \\[ x^{%d} y^{%d} \\]\n  \\eq\n  ", J1, J2);
    if (K[2]) {
      pflag = 1;  printf("4 b_1 \\beta_1^{%d}", J1+J2);
    }
    if (K[3]) {
      if (pflag) printf(" +\n  ");  pflag = 1;
      if (K[3] > 1) {
        printf("4 \\sum_{i=1}^{%d} c_i \\( \\gamma_i^{%d} \\delta_i^{%d} ", K[3], J1, J2);
        printf("+ \\gamma_i^{%d} \\delta_i^{%d} \\)", J2, J1);
      }
      else {
        printf("4 c_1 \\( \\gamma_1^{%d} \\delta_1^{%d} ", J1, J2);
        printf("+ \\gamma_1^{%d} \\delta_1^{%d} \\)", J2, J1);
      }
    }
    if (K[4]) {
      if (pflag) printf(" +\n  ");  pflag = 1;
      printf("8 d_1 \\epsilon_1^{%d}", J1+J2);
    }
    if (K[5]) {
      if (pflag) printf(" +\n  ");  pflag = 1;
      if (K[5] > 1) {
        printf("8 \\sum_{i=1}^{%d} e_i \\( \\zeta_i^{%d} + ", K[5], J1+J2);
        printf("\\zeta_i^{%d} \\eta_i^{%d} + \\zeta_i^{%d} \\eta_i^{%d} \\)", J1, J2, J2, J1);
      }
      else {
        printf("8 e_1 \\( \\zeta_1^{%d} + ", J1+J2);
        printf("\\zeta_1^{%d} \\eta_1^{%d} + \\zeta_1^{%d} \\eta_1^{%d} \\)", J1, J2, J2, J1);
      }
    }
    if (K[6]) {
      if (pflag) printf(" +\n  ");  pflag = 1;
      if (K[6] > 1) {
        printf("  \\\\\n  & & \\qquad8 \\sum_{i=1}^{%d} f_i \\( ", K[6]);
        printf("\\theta_i^{%d} \\mu_i^{%d} + ", J1, J2);
        printf("\\theta_i^{%d} \\mu_i^{%d} + ", J2, J1);
        printf("\\theta_i^{%d} \\lambda_i^{%d} + ", J1, J2);
        printf("\\theta_i^{%d} \\lambda_i^{%d} +", J2, J1);
        printf("\\mu_i^{%d} \\lambda_i^{%d} + ", J1, J2);
        printf("\\mu_i^{%d} \\lambda_i^{%d} \\)", J2, J1);
      }
      else {
        printf("  \\\\\n  & & \\qquad8 f_1 \\( ");
        printf("\\theta_1^{%d} \\mu_1^{%d} + ", J1, J2);
        printf("\\theta_1^{%d} \\mu_1^{%d} + ", J2, J1);
        printf("\\theta_1^{%d} \\lambda_1^{%d} + ", J1, J2);
        printf("\\theta_1^{%d} \\lambda_1^{%d} + ", J2, J1);
        printf("\\mu_1^{%d} \\lambda_1^{%d} + ", J1, J2);
        printf("\\mu_1^{%d} \\lambda_1^{%d} \\)", J2, J1);
      }
    }
    printf("\n  \\\\\n");
  } }

/* Deal with Subsystem III */
  if (m >= 3) {
  printf("\%% Subsystem III:\n");
  for (j1 = 1; j1 <= m; j1++)
  for (j2 = j1; j2 <= m-j1; j2++)
  for (j3 = j2; j3 <= m-j1-j2; j3++) {
    pflag = 0;  J1 = 2*j1;  J2 = 2*j2;  J3 = 2*j3;  eqno++;
    printf("  I \\[ x^{%d} y^{%d} z^{%d} \\]\n  \\eq\n  ", J1, J2, J3);
    if (K[4]) {
      pflag = 1;  printf("8 d_1 \\epsilon_1^{%d}", J1+J2+J3);
    }
    if (K[5]) {
      if (pflag) printf(" +\n  ");  pflag = 1;
      if (K[5] > 1) {
        printf("8 \\sum_{i=1}^{%d} e_i \\( \\zeta_i^{%d} \\eta_i^{%d} + ", K[5], J1+J2, J3);
        printf("\\zeta_i^{%d} \\eta_i^{%d} + ", J1+J3, J2);
        printf("\\zeta_i^{%d} \\eta_i^{%d} \\)", J2+J3, J1);
      }
      else {
        printf("8 e_1 \\( \\zeta_1^{%d} \\eta_1^{%d} + ", J1+J2, J3);
        printf("\\zeta_1^{%d} \\eta_1^{%d} + ", J1+J3, J2);
        printf("\\zeta_1^{%d} \\eta_1^{%d} \\)", J2+J3, J1);
      }
    }
    if (K[6]) {
      if (pflag) printf(" +\n  ");  pflag = 1;
      if (K[6] > 1) {
        printf("  \\\\\n  & & \\qquad8 \\sum_{i=1}^{%d} f_i \\( ", K[6]);
        printf("\\theta_i^{%d} \\mu_i^{%d} \\lambda_i^{%d} + ", J1, J2, J3);
        printf("\\theta_i^{%d} \\mu_i^{%d} \\lambda_i^{%d} + ", J1, J3, J2);
        printf("\\theta_i^{%d} \\mu_i^{%d} \\lambda_i^{%d} + ", J2, J1, J3);
        printf("\\theta_i^{%d} \\mu_i^{%d} \\lambda_i^{%d} + ", J2, J3, J1);
        printf("\\theta_i^{%d} \\mu_i^{%d} \\lambda_i^{%d} + ", J3, J1, J2);
        printf("\\theta_i^{%d} \\mu_i^{%d} \\lambda_i^{%d} \\)", J3, J2, J1);
      }
      else {
        printf("  \\\\\n  & & \\qquad8 f_1 \\( ");
        printf("\\theta_1^{%d} \\mu_1^{%d} \\lambda_1^{%d} + ", J1, J2, J3);
        printf("\\theta_1^{%d} \\mu_1^{%d} \\lambda_1^{%d} + ", J1, J3, J2);
        printf("\\theta_1^{%d} \\mu_1^{%d} \\lambda_1^{%d} + ", J2, J1, J3);
        printf("\\theta_1^{%d} \\mu_1^{%d} \\lambda_1^{%d} + ", J2, J3, J1);
        printf("\\theta_1^{%d} \\mu_1^{%d} \\lambda_1^{%d} + ", J3, J1, J2);
        printf("\\theta_1^{%d} \\mu_1^{%d} \\lambda_1^{%d} \\)", J3, J2, J1);
      }
    }
    printf("\n  \\\\\n");
  } }

  nflag = 0;
  printf("\%% Also Sprach:\n");
  if (K[1]) {
    printf("  1\n  \\eq\n  \\alpha_1\n");  nflag = 1;  eqno++;
  }
  if (K[2]) {
    if (nflag) printf("  \\\\\n");  nflag = 1;
    printf("  1/\\sqrt{2}\n  \\eq\n  \\beta_1\n");  eqno++;
  }
  if (K[3]) {
    if (nflag) printf("  \\\\\n");  nflag = 1;
    if (K[3] > 1)
      printf("  1\n  \\eq\n  \\gamma_i^2 + \\delta_i^2 \\qquad i = 1, \\dots, %d\n", K[3]);
    else
      printf("  1\n  \\eq\n  \\gamma_1^2 + \\delta_1^2\n");
    eqno += K[3];
  }
  if (K[4]) {
    if (nflag) printf("  \\\\\n");  nflag = 1;
    printf("  1/\\sqrt{3}\n  \\eq\n  \\epsilon_1\n");  eqno++;
  }
  if (K[5]) {
    if (nflag) printf("  \\\\\n");  nflag = 1;
    if (K[5] > 1)
      printf("  1\n  \\eq\n  2 \\zeta_i^2 + \\eta_i^2 \\qquad i = 1, \\dots, %d\n", K[5]);
    else
      printf("  1\n  \\eq\n  2 \\zeta_1^2 + \\eta_1^2\n");
    eqno += K[5];
  }
  if (K[6]) {
    if (nflag) printf("  \\\\\n");  nflag = 1;
    if (K[6] > 1) {
      printf("  1\n  \\eq\n  \\theta_i^2 + \\mu_i^2 + \\lambda_i^2 ");
      printf("\\qquad i = 1, \\dots, %d\n", K[6]);
    }
    else
      printf("  1\n  \\eq\n  \\theta_1^2 + \\mu_1^2 + \\lambda_1^2\n");
    eqno += K[6];
  }
  printf("\\ee\n\nThere are a total of %d equations.\n\n", eqno);
}
\end{verbatim} 
\normalsize\rm

%% file: C/eqnout.tex
Output from running `writestar': m = 4, K = 1 0 1 1 0 0:

\begin{eqnarray*}
  I [ 1 ]
  & = &
  6 a_1 +
  24 c_1 +
  8 d_1
  \\
  I [ x^{2} ]
  & = &
  2 a_1 \alpha_1^{2} +
  8 c_1 ( \gamma_1^{2} + \delta_1^{2} ) +
  8 d_1 \epsilon_1^{2}
  \\
  I [ x^{4} ]
  & = &
  2 a_1 \alpha_1^{4} +
  8 c_1 ( \gamma_1^{4} + \delta_1^{4} ) +
  8 d_1 \epsilon_1^{4}
  \\
  I [ x^{6} ]
  & = &
  2 a_1 \alpha_1^{6} +
  8 c_1 ( \gamma_1^{6} + \delta_1^{6} ) +
  8 d_1 \epsilon_1^{6}
  \\
  I [ x^{8} ]
  & = &
  2 a_1 \alpha_1^{8} +
  8 c_1 ( \gamma_1^{8} + \delta_1^{8} ) +
  8 d_1 \epsilon_1^{8}
  \\
  I [ x^{2} y^{2} ]
  & = &
  4 c_1 ( \gamma_1^{2} \delta_1^{2} + \gamma_1^{2} \delta_1^{2} ) +
  8 d_1 \epsilon_1^{4}
  \\
  I [ x^{2} y^{4} ]
  & = &
  4 c_1 ( \gamma_1^{2} \delta_1^{4} + \gamma_1^{4} \delta_1^{2} ) +
  8 d_1 \epsilon_1^{6}
  \\
  I [ x^{2} y^{6} ]
  & = &
  4 c_1 ( \gamma_1^{2} \delta_1^{6} + \gamma_1^{6} \delta_1^{2} ) +
  8 d_1 \epsilon_1^{8}
  \\
  I [ x^{4} y^{4} ]
  & = &
  4 c_1 ( \gamma_1^{4} \delta_1^{4} + \gamma_1^{4} \delta_1^{4} ) +
  8 d_1 \epsilon_1^{8}
  \\
  I [ x^{2} y^{2} z^{2} ]
  & = &
  8 d_1 \epsilon_1^{6}
  \\
  I [ x^{2} y^{2} z^{4} ]
  & = &
  8 d_1 \epsilon_1^{8}
  \\
  1
  & = &
  \alpha_1
  \\
  1
  & = &
  \gamma_1^2 + \delta_1^2
  \\
  1/\sqrt{3}
  & = &
  \epsilon_1
\end{eqnarray*}

There are a total of 14 equations.

%% file: matlab/cubature.tex
\small\tt 
\begin{verbatim}
% Solve the system (*) of moment equations, for input parameters
% m and K, to obtain a cubature rule of degree 2m+1.
%
%  David  De Wit
%  September 20  --  October 27  1993

function x = cubature(m, j, t)

format compact;      format long

if (m == 1)
  K = [1 0 0 0 0 0];
elseif (m == 2)
  K = [1 0 0 1 0 0];
elseif (m == 3)
  K = [1 1 0 1 0 0];
elseif (m == 4)
% K = [1 0 0 1 1 0];
  K = [1 0 1 1 0 0];
elseif (m == 5)
  K = [1 1 0 1 1 0];
elseif (m == 6)
  K = [1 1 1 1 1 0];
% K = [1 0 1 0 2 0];
elseif (m == 7)
  K = [1 0 1 1 2 0];
elseif (m == 8)
  K = [1 0 1 1 3 0];
elseif (m == 9)
  K = [1 1 0 1 3 1];
% K = [1 1 1 1 2 1];
elseif (m == 10)
  K = [1 1 1 1 3 1];
% K = [1 1 2 1 2 1];
end

M = [2 2 3 2 3 4]';
L(1) = 0;
for i = 2:6
  L(i) = L(i-1) + M(i-1)*K(i-1);
end

x0 = rand(K*M,1);
if (K(1)), x0(L(1)+2) = 1; end
if (K(2)), x0(L(2)+2) = 1/sqrt(2); end
if (K(4)), x0(L(4)+2) = 1/sqrt(3); end

options = [1 eps eps 0 1];  options(14) = 10000;

%  Pass the values to the solution routine.
x = fsolve('momenteq', x0, options, [], m, K, L);

disp('m, j, t, K are:'); disp([m j t K]);

%  Decode the resulting x into known variables:
if (K(1))
  a = x(1), alpha = x(2)
end
if (K(2))
 b = x(L(2)+1), beta = x(L(2)+2)
end
if (K(3))
  c       = x(L(3)+1:L(3)+K(3))
  gamma   = x(L(3)+K(3)+1:L(3)+2*K(3));    delta = x(L(3)+2*K(3)+1:L(3)+3*K(3))
end
if (K(4))
  d = x(L(4)+1), epsilon = x(L(4)+2)
end
if (K(5))
  e       = x(L(5)+1:L(5)+K(5))
  zeta    = x(L(5)+K(5)+1:L(5)+2*K(5));      eta = x(L(5)+2*K(5)+1:L(5)+3*K(5))
end
if (K(6))
  f       = x(L(6)+1:L(6)+K(6));           theta = x(L(6)+K(6)+1:L(6)+2*K(6))
  mu      = x(L(6)+2*K(6)+1:L(6)+3*K(6)); lambda = x(L(6)+3*K(6)+1:L(6)+4*K(6))
end
\end{verbatim} 
\normalsize\rm

%% file: matlab/momenteq.tex
\small\tt 
\begin{verbatim}
/*
  C code that compiles into a MATLAB .mex file that `evaluates'
  the moment equations (*).

  David  De Wit and Martin Sharry and Peter Adams
  September 24  --  September 29  1993
*/

#include <math.h>
#include "mex.h"

/* Input and Output Arguments */

#define  x_IN  prhs[0]
#define  m_IN  prhs[1]
#define  K_IN  prhs[2]
#define  L_IN  prhs[3]
#define  F_OUT  plhs[0]

static int Neq[10] = {2, 4, 7, 11, 16, 23, 31, 41, 53, 67};
static int M[6] = {2, 2, 3, 2, 3, 4};

mexFunction(nlhs, plhs, nrhs, prhs)
int    nlhs, nrhs;
Matrix *plhs[], *prhs[];
{
  double *F, *x, m, *K, *L;
  int    i, mm; 

  /* Check for proper number of arguments */
  if (nrhs != 4) {
    mexErrMsgTxt("momenteq requires four input arguments.");
  } else if (nlhs > 1) {
    mexErrMsgTxt("momenteq requires one output argument.");
  }

  /* Assign pointers to the various parameters */

  x = mxGetPr(x_IN);  m = mxGetScalar(m_IN);
  K = mxGetPr(K_IN);  L = mxGetPr(L_IN);

  mm = Neq[(int) m - 1];
  for (i = 0; i < 6; i++)
    mm += K[i]*(M[i] - 1);
  F_OUT = mxCreateFull(mm, 1, REAL); 
  F = mxGetPr(F_OUT);

  momenteq(F, x, m, K, L);
}

#define  pi  3.141592653589793238462643383280
double   mom[4][6][11];

void setmom()
{
  mom[0][0][0]  =  4*pi;            mom[0][0][1]  =  4*pi/3;
  mom[0][0][2]  =  4*pi/5;          mom[0][0][3]  =  4*pi/7;
  mom[0][0][4]  =  4*pi/9;          mom[0][0][5]  =  4*pi/11;
  mom[0][0][6]  =  4*pi/13;         mom[0][0][7]  =  4*pi/15;
  mom[0][0][8]  =  4*pi/17;         mom[0][0][9]  =  4*pi/19;
  mom[0][0][10] =  4*pi/21;         mom[0][1][1]  =  4*pi/15;
  mom[0][1][2]  =  4*pi/35;         mom[0][1][3]  =  4*pi/63;
  mom[0][1][4]  =  4*pi/99;         mom[0][1][5]  =  4*pi/143;
  mom[0][1][6]  =  4*pi/195;        mom[0][1][7]  =  4*pi/255;
  mom[0][1][8]  =  4*pi/323;        mom[0][1][9]  =  4*pi/399;
  mom[0][2][2]  =  4*pi/105;        mom[0][2][3]  =  4*pi/231;
  mom[0][2][4]  =  4*pi/429;        mom[0][2][5]  =  4*pi/715;
  mom[0][2][6]  =  4*pi/1105;       mom[0][2][7]  =  4*pi/1615;
  mom[0][2][8]  =  4*pi/2261;       mom[0][3][3]  =  20*pi/3003;
  mom[0][3][4]  =  4*pi/1287;       mom[0][3][5]  =  4*pi/2431;
  mom[0][3][6]  =  4*pi/4199;       mom[0][3][7]  =  4*pi/6783;
  mom[0][4][4]  =  28*pi/21879;     mom[0][4][5]  =  28*pi/46189;
  mom[0][4][6]  =  4*pi/12597;      mom[0][5][5]  =  12*pi/46189;
  mom[1][1][1]  =  4*pi/105;        mom[1][1][2]  =  4*pi/315;
  mom[1][1][3]  =  4*pi/693;        mom[1][1][4]  =  4*pi/1287;
  mom[1][1][5]  =  4*pi/2145;       mom[1][1][6]  =  4*pi/3315;
  mom[1][1][7]  =  4*pi/4845;       mom[1][1][8]  =  4*pi/6783;
  mom[1][2][2]  =  4*pi/1155;       mom[1][2][3]  =  4*pi/3003;
  mom[1][2][4]  =  4*pi/6435;       mom[1][2][5]  =  4*pi/12155;
  mom[1][2][6]  =  4*pi/20995;      mom[1][2][7]  =  4*pi/33915;
  mom[1][3][3]  =  4*pi/9009;       mom[1][3][4]  =  4*pi/21879;
  mom[1][3][5]  =  4*pi/46189;      mom[1][3][6]  =  4*pi/88179;
  mom[1][4][4]  =  28*pi/415701;    mom[1][4][5]  =  4*pi/138567;
  mom[2][2][2]  =  4*pi/5005;       mom[2][2][3]  =  4*pi/15015;
  mom[2][2][4]  =  4*pi/36465;      mom[2][2][5]  =  12*pi/230945;
  mom[2][2][6]  =  4*pi/146965;     mom[2][3][3]  =  4*pi/51051;
  mom[2][3][4]  =  4*pi/138567;     mom[2][3][5]  =  4*pi/323323;
  mom[2][4][4]  =  4*pi/415701;     mom[3][3][3]  =  20*pi/969969;
  mom[3][3][4]  =  20*pi/2909907;
}

momenteq(F, x, m, K, L)
double  F[], x[], m, K[], L[];
{
  int        iK[6], iL[6], i, im, j1, j2, j3, p;
  double     *a, *alpha, *b, *beta, *c,  *gammah, *delta, *d, *epsilon,
             *e, *zeta, *eta, *f, *theta, *mu, *lambda, J1, J2, J3;
  static int havesetmom = 0;

  if(!havesetmom) {
    setmom();  havesetmom = 1;
  }
  im = (int) m;
  for (i = 0; i < 6; i++) {
    iK[i] = (int) K[i];  iL[i] = (int) L[i];
  }

  a       = (double *) malloc ( iK[0] * sizeof(double) );
  alpha   = (double *) malloc ( iK[0] * sizeof(double) );
  b       = (double *) malloc ( iK[1] * sizeof(double) );
  beta    = (double *) malloc ( iK[1] * sizeof(double) );
  c       = (double *) malloc ( iK[2] * sizeof(double) );
  gammah  = (double *) malloc ( iK[2] * sizeof(double) );
  delta   = (double *) malloc ( iK[2] * sizeof(double) );
  d       = (double *) malloc ( iK[3] * sizeof(double) );
  epsilon = (double *) malloc ( iK[3] * sizeof(double) );
  e       = (double *) malloc ( iK[4] * sizeof(double) );
  zeta    = (double *) malloc ( iK[4] * sizeof(double) );
  eta     = (double *) malloc ( iK[4] * sizeof(double) );
  f       = (double *) malloc ( iK[5] * sizeof(double) );
  theta   = (double *) malloc ( iK[5] * sizeof(double) );
  mu      = (double *) malloc ( iK[5] * sizeof(double) );
  lambda  = (double *) malloc ( iK[5] * sizeof(double) );

  for (i = 0; i < iK[0]; i++) {
    a[i] = x[iL[0]+i];   alpha[i] = x[iL[0]+1*iK[0]+i];
  }
  for (i = 0; i < iK[1]; i++) {
   b[i] = x[iL[1]+i];     beta[i] = x[iL[1]+1*iK[1]+i];
  }
  for (i = 0; i < iK[2]; i++) {
    c[i] = x[iL[2]+i];  gammah[i] = x[iL[2]+1*iK[2]+i]; delta[i] = x[iL[2]+2*iK[2]+i];
  }
  for (i = 0; i < iK[3]; i++) {
    d[i] = x[iL[3]+i]; epsilon[i] = x[iL[3]+1*iK[3]+i];
  }
  for (i = 0; i < iK[4]; i++) {
    e[i] = x[iL[4]+i];    zeta[i] = x[iL[4]+1*iK[4]+i]; eta[i] = x[iL[4]+2*iK[4]+i];
  }
  for (i = 0; i < iK[5]; i++) {
    f[i] = x[iL[5]+i];   theta[i] = x[iL[5]+1*iK[5]+i]; mu[i] = x[iL[5]+2*iK[5]+i];
                        lambda[i] = x[iL[5]+3*iK[5]+i];
  }

/* Subsystem I */
  F[p = 0] = - mom[0][0][0];
  for (i = 0; i < iK[0]; i++) F[p] +=  6*a[i];
  for (i = 0; i < iK[1]; i++) F[p] += 12*b[i];
  for (i = 0; i < iK[2]; i++) F[p] += 24*c[i];
  for (i = 0; i < iK[3]; i++) F[p] +=  8*d[i];
  for (i = 0; i < iK[4]; i++) F[p] += 24*e[i];
  for (i = 0; i < iK[5]; i++) F[p] += 48*f[i];

  for (j1 = 1; j1 <= im; j1++) {
    J1 = 2.0*j1;
    F[++p] = - mom[0][0][j1];
    for (i = 0; i < iK[0]; i++)
      F[p] +=  2*a[i]*pow(alpha[i],J1);
    for (i = 0; i < iK[1]; i++)
      F[p] +=  8*b[i]*pow(beta[i],J1);
    for (i = 0; i < iK[2]; i++)
      F[p] +=  8*c[i]*(pow(gammah[i],J1) + pow(delta[i],J1));
    for (i = 0; i < iK[3]; i++)
      F[p] +=  8*d[i]*pow(epsilon[i],J1);
    for (i = 0; i < iK[4]; i++)
      F[p] +=  8*e[i]*(2*pow(zeta[i],J1) + pow(eta[i],J1));
    for (i = 0; i < iK[5]; i++)
      F[p] += 16*f[i]*(pow(theta[i],J1) + pow(mu[i],J1) + pow(lambda[i],J1));
  }

/* Subsystem II  */
  if (im >= 2)
    for (j1 = 1; j1 <= im; j1++) {
      J1 = 2.0*j1;
      for (j2 = j1; j2 <= im-j1; j2++) {
        J2 = 2.0*j2;
        F[++p] = - mom[0][j1][j2];
        for (i = 0; i < iK[1]; i++)
          F[p] +=  4*b[i]*pow(beta[i],J1+J2);
        for (i = 0; i < iK[2]; i++)
          F[p] +=  4*c[i]*(
                   pow(gammah[i],J1)*pow(delta[i],J2) + pow(gammah[i],J2)*pow(delta[i],J1));
        for (i = 0; i < iK[3]; i++)
          F[p] +=  8*d[i]*pow(epsilon[i],J1+J2);
        for (i = 0; i < iK[4]; i++)
          F[p] +=  8*e[i]*(
                   pow(zeta[i],J1+J2) +
                   pow(zeta[i],J1)*pow(eta[i],J2) +
                   pow(zeta[i],J2)*pow(eta[i],J1));
        for (i = 0; i < iK[5]; i++)
          F[p] += 8*f[i]*(
                  pow(theta[i],J1)*pow(mu[i],J2) + pow(theta[i],J2)*pow(mu[i],J1) +
                  pow(theta[i],J1)*pow(lambda[i],J2) + pow(theta[i],J2)*pow(lambda[i],J1) +
                  pow(mu[i],J1)*pow(lambda[i],J2) + pow(mu[i],J2)*pow(lambda[i],J1));
    } }

/* Subsystem III  */
  if (im >= 3)
    for (j1 = 1; j1 <= im; j1++) {
      J1 = 2.0*j1;
      for (j2 = j1; j2 <= im-j1; j2++) {
        J2 = 2.0*j2;
        for (j3 = j2; j3 <= im-j1-j2; j3++) {
          J3 = 2.0*j3;
          p++;
          F[p] = - mom[j1][j2][j3];
          for (i = 0; i < iK[3]; i++)
            F[p] +=  8*d[i]*pow(epsilon[i],J1+J2+J3);
          for (i = 0; i < iK[4]; i++)
            F[p] +=  8*e[i]*(
                     pow(zeta[i],J1+J2)*pow(eta[i],J3) +
                     pow(zeta[i],J1+J3)*pow(eta[i],J2) +
                     pow(zeta[i],J2+J3)*pow(eta[i],J1));
          for (i = 0; i < iK[5]; i++)
            F[p] += 8*f[i]*(
                    pow(theta[i],J1)*pow(mu[i],J2)*pow(lambda[i],J3) +
                    pow(theta[i],J1)*pow(mu[i],J3)*pow(lambda[i],J2) +
                    pow(theta[i],J2)*pow(mu[i],J1)*pow(lambda[i],J3) +
                    pow(theta[i],J2)*pow(mu[i],J3)*pow(lambda[i],J1) +
                    pow(theta[i],J3)*pow(mu[i],J1)*pow(lambda[i],J2) +
                    pow(theta[i],J3)*pow(mu[i],J2)*pow(lambda[i],J1));
    } } };

/* Also Sprach  */
  for (i = 0; i < iK[0]; i++)
    F[++p] = alpha[i] - 1;
  for (i = 0; i < iK[1]; i++)
    F[++p] = beta[i] - 1/sqrt(2.0);
  for (i = 0; i < iK[2]; i++)
    F[++p] = pow(gammah[i],2.0) + pow(delta[i],2.0) - 1;
  for (i = 0; i < iK[3]; i++)
    F[++p] = epsilon[i] - 1/sqrt(3.0);
  for (i = 0; i < iK[4]; i++)
    F[++p] = 2*pow(zeta[i],2.0) + pow(eta[i],2.0) - 1;
  for (i = 0; i < iK[5]; i++)
    F[++p] = pow(theta[i],2.0) + pow(mu[i],2.0) + pow(lambda[i],2.0) - 1;

  free(a); free(alpha); free(b); free(beta); free(c); free(gammah); free(delta);
  free(d); free(epsilon); free(e); free(zeta); free(eta);
  free(f); free(theta); free(mu); free(lambda);
}
\end{verbatim} 
\normalsize\rm

%% file: matlab/u3prod.tex
\small\tt 
\begin{verbatim}
function [y, A] = u3prod(M)

% function [y, A] = u3prod(M)
%
% Product rule cubature for U_3. On input, M is the number of points in
% the basic rules. Gives a Gauss product rule of degree 2M - 1 on 2 M^2
% Cartesian points -+ y, with weights A. y is an M^2 x 3 matrix, A an
% M^2 column vector.  To obtain a rule of degree p, a rule on
% (p + 1)^2 / 2 points is required.  Requires gauss.m, which returns
% the points and weights of the $ M $-point 1D Gauss-Legendre rule.
% See Stroud (1971), p 41.
%
% David  De Wit
% February 11  -  February 12  1993

if ~exist('M'), M = 5; end

[yG, AG] = gauss(M,-1,1);
j = [1:M]';

yCI = cos((2*j - 1)*pi/(2*M));    ACI = ones(size(j))*pi/M;

y1 = sqrt(1 - yCI.^2) * sqrt(1 - yG.^2)';
y2 = yCI * sqrt(1 - yG.^2)';
y3 = ones(M,1) * yG';
A = ACI * AG';

y = [y1(:) y2(:) y3(:)];          A = A(:);

plot3(y(:,1),y(:,2),y(:,3),'+r');

% Test by integrating all functions of x, y and z of degree < M.
% Works beautifully!

% t = M-1;        m = 1;
% for j = 0:t, for k = j:t, for l = k:t
%   if (j+k+l <= t)
%     i = [j k l];
%     tabi(m,:) = i;
%     table(m) = 2*(y(:,1).^(2*i(1)).*y(:,2).^(2*i(2)) ...
%                .*y(:,3).^(2*i(3)))'*A - ...
%                2*prod(gamma(i+1/2))/gamma(3/2+sum(i));
%     m = m + 1;
%   end
% end, end, end
% tabi
% table'
% m
\end{verbatim} 
\normalsize\rm

%% file: matlab/gauss.tex
\small\tt 
\begin{verbatim}
function [x, w] = gauss(n, a, b)

% function [x, w] = gauss(n, a, b)
% Returns {x, w}, the weights and of the n-point Gauss- Legendre
% quadrature rule on the interval  [a, b]
%
% Graeme Chandler 1992

if (n==1)
    x = (a+b)/2 ; w = b-a ; return
end

m = 1:2:2*n-1 ;
m = (1:n-1) ./ sqrt(m(1:n-1).*m(2:n));   %  m is off-diagonal of matrix
[w, x] = eig(diag(m,-1)+diag(m,1));      %  Find spectrum of matrix
x = (a+b)/2 + ((b-a)/2)*diag(x);         %  x are the eigenvalues
w = (b-a)*(w(1,:).^2)';                  %  w from first components
[x, m] = sort(x);        w = w(m);       %  ascending order
\end{verbatim} 
\normalsize\rm

%% file: Mathematica/IPPBasicSoln.tex
\small\tt 
\begin{verbatim}
(*
  Find basic solutions to the IPP, ignoring the integer requirement,
  for m = 1, ..., 20.

  David  De Wit
  August 4  1993
*)

mmax = 20
A = {
  { 1,  1,  2,  1,  2,  3},
  { 0,  0,  0,  1,  2,  3},
  { 0,  0,  2,  0,  0,  3},
  { 0,  0,  0,  0,  0,  3},
  {-1,  0,  0,  0,  0,  0},
  { 0, -1,  0,  0,  0,  0},
  { 0,  0,  0, -1,  0,  0}
}
c = {6, 12, 24, 8, 24, 48}
b = Transpose[ {
  { 0,  1,  2,  3,  4,  5,  7,  8, 10, 12, 14, 16, 19, 21, 24, 27, 30, 33, 37, 40, 44},
  { 0,  0,  0,  1,  1,  2,  3,  4,  5,  7,  8, 10, 12, 14, 16, 19, 21, 24, 27, 30, 33},
  { 0,  0,  0,  0,  0,  0,  1,  1,  2,  3,  4,  5,  7,  8, 10, 12, 14, 16, 19, 21, 24},
  { 0,  0,  0,  0,  0,  0,  0,  0,  0,  1,  1,  2,  3,  4,  5,  7,  8, 10, 12, 14, 16},
  {-1, -1, -1, -1, -1, -1, -1, -1, -1, -1, -1, -1, -1, -1, -1, -1, -1, -1, -1, -1, -1},
  {-1, -1, -1, -1, -1, -1, -1, -1, -1, -1, -1, -1, -1, -1, -1, -1, -1, -1, -1, -1, -1},
  {-1, -1, -1, -1, -1, -1, -1, -1, -1, -1, -1, -1, -1, -1, -1, -1, -1, -1, -1, -1, -1}
}]
sol = Table[0, {i, mmax+1}, {j, 6}]
nmin = Table[0, {i, mmax+1}]
For [m = 1, m <= mmax+1, m++,
  sol[[m]] = LinearProgramming[c, A, b[[m]]]
]
nmin = c . Transpose[sol]

stmp = OpenWrite["IPPBasicSolnData.tex.bak"]
For [m = 1, m <= mmax+1, m++,
  WriteString[{"stdout", stmp}, "$ ", m-1, " $ & $ ", nmin[[m]]
  ];
  For [n = 1, n <=6, n++,
    WriteString [{"stdout", stmp}, " $ & $ ", InputForm[sol[[m,n]]]
    ];
  ];
  WriteString[ {"stdout", stmp}, " $ \\\\\n" ];
];
Close[stmp]
\end{verbatim} 
\normalsize\rm

%% file: Mathematica/U3Moments.tex
\small\tt 
\begin{verbatim}
(*
  Find the moments for integrals of polynomials over U_3, and output
  into a file suitable for inclusion as a LaTeX table. The output is
  also included in the programs momenteq.c and writestar.c. The method
  is extremely crude, but operational!

  David  De Wit
  August 4  -- August 12  1993
*)

stmp = OpenWrite["U3MomentData.tex"]
For [j1 = 0, j1 <= 10, j1++,
For [j2 = j1, j2 <= 10-j1, j2++,
For [j3 = j2, j3 <= 10-j2-j1, j3++,
  f[t_] = (Cos[t])^(2 j1) (Sin[t])^(2 j2);
  g[t_] = (Cos[t])^(2 j1 + 2 j2 + 1) (Sin[t])^(2 j3);
  Int = Integrate[g[t], {t,-Pi/2,Pi/2}] Integrate[f[t], {t,-Pi,Pi}];
  WriteString[
    {"stdout", stmp},
    "$ ", j1, " $ & $ ", j2, " $ & $ ", j3, " $ & $ ", InputForm[Int], " $ \\\\", "\n"]
] ] ]
Close[stmp]
\end{verbatim} 
\normalsize\rm